\documentclass[11pt]{article}
\usepackage[ngerman,english]{babel}
\usepackage[utf8]{inputenc}
\usepackage{amssymb,amsmath,amstext,amsthm,amsfonts, ams fonts}\usepackage{colortbl,color}
\usepackage{graphicx,url}
\usepackage{hyperref}
\usepackage[font=small,skip=0pt]{caption}
\usepackage[authoryear,sort&compress]{natbib}
\usepackage{booktabs}
\usepackage[flushmargin,hang]{footmisc}
\usepackage{setspace}
\usepackage[tight,small]{subfigure}
\usepackage[verbose=true]{microtype}
\usepackage{multirow,tabularx,rotating,multicol}
\usepackage{times,dsfont}
\usepackage{bbm}
\usepackage{nicefrac}
\usepackage{lineno}
\usepackage{calc}
\usepackage{changes}
\graphicspath{{pics/}}

\definecolor{cWhite}{rgb}{1,1,1}
\definecolor{cLightRed}{rgb}{1,.70,.70}
\definecolor{cLightYellow}{rgb}{.90,.85,.55}
\definecolor{cLightGray}{rgb}{.90,.90,.90}
\definecolor{cMediumGray}{rgb}{.70,.70,.70}
\definecolor{cDarkGray}{rgb}{.50,.60,.70}
\definecolor{darkblue}{rgb}{0,0,.6}
\hypersetup{%
pdftitle={}, 
pdfproducer={PDFTeX}, 
pdfkeywords={Working Paper}, 
plainpages=false,
pdfstartpage={1}, 
pdfpagelayout=OneColumn, 
pdfpagemode=UseOutlines, 
linkcolor=darkblue, 
citecolor=darkblue,
urlcolor=black, 
breaklinks=true, 
colorlinks=true, 
citebordercolor=0 0 0, 
filebordercolor=0 0 0, 
linkbordercolor=0 0 0, 
menubordercolor=0 0 0, %
urlbordercolor=0 0 0, 
pdfhighlight=/I, %
pdfborder=0 0 0, 
bookmarksopen=true, %
bookmarksnumbered=true %
}
\allowdisplaybreaks

\renewcommand{\P}{\mathbb{P}}

\newcommand{\E}{\mathbb{E}}

\newcommand{\R}{\mathds{R}}
\newcommand{\1}{\mathbbm{1}}

\newcommand{\KLEINO}{{\scriptstyle{\mathcal{O}}}}
\DeclareMathAccent{\verywidehat}{\mathord}{largesymbols}{'144}

\newcommand{\Cov}{\mathbb{C}\textnormal{O\hspace*{0.02cm}V}}
\newcommand{\cov}{\mathbb{C}\textnormal{o\hspace*{0.02cm}v}}

\renewcommand{\:}{\mathrel{\mathop{:}}}

\newtheorem{theo}{Theorem}
\newtheorem{corro}{Corollary}

\newtheorem{sassump}{Assumption}

 \newcommand{\defeq}{\mathrel{\mathop:}=}

\newcommand{\comment}[1]{ }

\usepackage{verbatim,color}
\def\boxit#1{\vbox{\hrule\hbox{\vrule\kern6pt
          \vbox{\kern6pt#1\kern6pt}\kern6pt\vrule}\hrule}}

\title{Estimating the Spot Covariation of Asset Prices -- Statistical Theory and Empirical Evidence\thanks{
Financial support from the Deutsche Forschungsgemeinschaft via SFB 649 ``Economic Risk'' and FOR 1735 ``Structural Inference in Statistics: Adaptation and Efficiency'' is gratefully acknowledged. Hautsch also acknowledges financial support from the Wiener Wissenschafts-, Forschungs- und Technologiefonds (WWTF). Malec thanks the Cambridge INET for financial support.}}
\author{{Markus Bibinger}\thanks{Department of Mathematics and Computer Science, University of Marburg. Email: bibinger@mathematik.uni-marburg.de. Address:  Hans-Meerwein-Stra\ss e~6, D-35032 Marburg, Germany.}\hspace{2em} {Nikolaus Hautsch}\thanks{Faculty of Business, Economics and Statistics,  University of Vienna and Center for Financial Studies, Frankfurt. Email: nikolaus.hautsch@univie.ac.at. Address: Oskar-Morgenstern-Platz~1, A-1090 Vienna, Austria.} \hspace{2em} {Peter Malec}\thanks{Faculty of Economics, University of Cambridge. Email: pm563@cam.ac.uk. Address: Sidgwick Avenue, Cambridge CB3 9DD, United Kingdom.} \hspace{2em} {Markus Reiss}\thanks{Institute of Mathematics, Humboldt-Universität zu Berlin. Email: mreiss@math.hu-berlin.de. Address: Unter den Linden 6, D-10099 Berlin, Germany.}}

\date{}
\begin{document}
\maketitle

\begin{abstract}
\noindent  We propose a new estimator for the spot covariance matrix of a multi-dimensional continuous semi-martingale log asset price process which is subject to noise and non-synchronous observations. The estimator is constructed based on a local average of block-wise parametric spectral covariance estimates. The latter originate from a local method of moments (LMM) which recently has been introduced by \citet{BHMR}. We prove consistency and a point-wise stable central limit theorem for the proposed spot covariance estimator in a very general setup with stochastic volatility, leverage effects and  general noise distributions. Moreover, we extend the LMM estimator to be robust against autocorrelated noise and propose a method to adaptively infer the autocorrelations from the data. Based on  simulations we provide empirical guidance on the effective implementation of the estimator and apply it to high-frequency data of a cross-section of Nasdaq blue chip stocks. Employing the estimator to estimate spot covariances, correlations and volatilities in normal but also unusual periods yields novel insights into intraday covariance and correlation dynamics. We show that intraday (co-)variations (i) follow underlying periodicity patterns, (ii) reveal substantial intraday variability associated with (co-)variation risk, and (iii) can increase strongly and nearly instantaneously if new information arrives.
%
\end{abstract}


\section{Introduction\label{sec:intro}}

Recent literature in financial econometrics and empirical finance reports strong empirical evidence for distinct time variations in daily and long-term correlations between asset prices. While the literature proposes several approaches to estimate spot \emph{variances}, 
there is a lack of empirical approaches and corresponding statistical theory to estimate spot \emph{covariances} using high-frequency data. In this paper, we aim at filling this gap in the literature and propose a novel estimator for the spot covariance matrix of a multi-dimensional continuous semi-martingale log asset price process which is observed at non-synchronous times under noise. 

Our study is mainly related to two fields of literature. First, there is a vast body of papers on the estimation of integrated covariance matrices, while accounting for market microstructure noise and the asynchronicity of observations. Starting from the seminal realized covariance estimator by \citet{barnd2004} which neglects both types of frictions, \citet{hy3} propose a consistent and efficient estimator under asynchronicity, but in the absence of microstructure noise. Estimators accounting for both types of frictions are, among others, the quasi maximum likelihood estimator by \citet{Ait-Sahalia2010}, the multivariate realized kernel estimator by \citet{bn1}, the multivariate pre-averaging estimator by \citet{PHY}, the two-scale estimator by \citet{lancov2011}, and the LMM estimator by \citet{BHMR}. \citet{aits2015} show how to estimate diagonalized integrated covariance matrices exploiting methods from \citet{jaco2013} to deal with functional transformations of volatility. 

Second, there is considerable literature on spot volatility estimation. A nonparametric (kernel-type) estimator in the absence of microstructure noise is put forward by \citet{fost96}, \citet{fan08} and \citet{kristensen}. To account for noise, the predominant approach is to compute a difference quotient based on a noise-robust integrated volatility estimator, e.g., the (univariate) realized kernel, the pre-averaging estimator or the two-scale estimator. Here, examples include \citet{mykl08}, \citet{manc15}, \citet{bos12} and \citet{zu14}. An alternative approach based on series estimators of \textit{non-stochastic} spot volatility  is introduced by \citet{munk2010b}, while \citet{munk2010} study optimal convergence rates in the aforementioned setting. Similarly, rates for the stochastic volatility case are derived by  \citet{hoff2012} who also propose a wavelet-type estimator attaining this rate. Finally, estimators that are robust to jumps, but neglect microstructure noise are put forward, e.g., by \citet{ait-sahalia09}, \citet{ander09} and \citet{band09}. Spot volatility estimation is relevant also for multiple-step approaches to perform inference on functionals that hinge on the volatility, see, e.g., \citet{kaln2016} and \citet{Li2016}. In the same way, our theory provides a foundation for multi-dimensional multiple-step inference based on the volatility matrix.

Interestingly, the problem of estimating the spot \textit{covariance matrix} in the presence of microstructure noise and asynchronicity effects has not yet been addressed in a study on its own. Our paper thus bridges the gap between the two fields of literature outlined above. In this context, note that spot covariance estimates are not derived as direct extensions of variance estimates under asynchronicity.
Our estimator is constructed based on local averages of block-wise parametric spectral covariance estimates. The latter are estimated employing the local method of moments (LMM) estimator proposed by \citet{BHMR}, which has been shown to attain the optimal rate and, moreover, a statistical lower bound for the asymptotic variance for the estimation of the integrated covariance matrix. As the LMM estimator builds on locally constant approximations of the underlying covariance process and estimates them block-wise, it provides a natural setting to construct a spot covariance estimator. Our methodological contribution is as follows: First, we construct the new spot covariance matrix estimator. Second, we derive a stable central limit theorem, showing the consistency and asymptotic normality of this estimator. For both, we consider a more realistic model than previous works based on the LMM method, for instance, allowing for autocorrelated market microstructure noise. 

Compared to integrated (co-)variance estimators, spot estimators inherently feature slower convergence rates due to the additional smoothing involved. We prove that our spot estimator can attain rate-optimality and satisfies a point-wise stable central limit theorem at almost optimal rate. Moreover, the reduced variance effect of the LMM estimator due to multivariate weight matrices carries over to spot covariance matrix estimation and appears to be relevant in practice.
Finally, as reported by \citet{hans06}, \citet{zhangmykland2}, and \citet{ACH2017} among others, microstructure noise appears to violate the traditional i.i.d.~assumption, exhibiting more complex dependence structures. Adjusting both our spot covariance estimator as well as the original LMM integrated covariance estimator by \cite{BHMR} to incorporate noise autocorrelation in a robust manner is an important extension, which makes the use of the methods in applications more attractive.

The approach presented here does not account for jumps in the log-price process. From a methodological point of view, an extension to disentangle jumps and continuous components utilizing a truncation technique as in \cite{bibwink2015} appears feasible. In the given framework, however, due to additional tuning parameters involved this would require a comprehensive extension, which would dilute the main new estimation ideas. Consequently, our proposed spot covariance estimator does not separate between a diffusive and jump component. For our empirical results and corresponding conclusions, this is not a limitation since in any case, potential jumps are consistently captured by the spot estimator. Moreover, \citet{chris14} show that, when considering  data sampled at the tick-by-tick level, jumps are detected far less frequently than based on a coarser sampling grid. 

Our approach allows for an efficient recovery of latent \emph{intraday spot (co-)volatility paths} of individual stocks. We provide simulation-based evidence on an effective implementation of the estimator depending on the choice of underlying smoothing parameters. Empirical studies on the role of high-frequency trading, the impact of market fragmentation and the usefulness of volatility circuit breakers might heavily benefit from the availability of high-frequency covariance estimators which are applicable in higher dimensions. Further, spot covariance estimates are a necessary building block for co-jump tests, see \citet{bibwink2015}. Finally, an important objective of this paper is to provide first empirical evidence on the intraday behavior of spot covariances and correlations.
Applying the spot covariance estimator to four years of quote data for 30 of the most liquid constituents of the Nasdaq 100, we obtain novel empirical findings. First, there is a distinct intraday seasonality pattern going beyond volatilities as covariances  exhibit a U-shape and correlations increase throughout the day. Second, spot (co-)variation reveals substantial intraday variability and thus reflect (co-)variation risk. Finally, spot covariances and correlations change substantially during flash crashes or the arrival of fundamental information. 

The remainder of the paper is structured as follows. Section \ref{sec:data} states a brief description of the data and empirical objectives. Section~\ref{sec:estimspotcov} theoretically introduces the proposed spot estimator and gives its asymptotic properties in detail. In Section~\ref{sec:sim}, we present a simulation study analyzing the estimator's sensitivity to the choice of input parameters and demonstrating its finite sample accuracy. Section~\ref{sec:emp} provides empirical evidence on spot (co-)variances, correlations and volatilities based on Nasdaq data. Finally, Section~\ref{sec:concl} concludes. Supplementary material is contained in a web appendix available on \url{https://www.mathematik.uni-marburg.de/~stochastik/material/Web_Appendix.pdf}.

\section{High-Frequency Data and Spot Covariances} \label{sec:data}
We employ ask and bid quote data for $30$ of the most liquid constituents of the Nasdaq 100 index. The sample period is from May 2010 to April 2014. The underlying data is provided by the LOBSTER database. The latter reconstructs the order book from a message stream, which is part of Nasdaq's historical TotalView-ITCH data and contains all limit order submissions, cancellations and executions on each trading day \citep[see][]{huang2011} on the Nasdaq market. Accordingly, the corresponding transaction data can be read out from the above message files directly. For the resulting datasets all recorded events are time stamped with at least millisecond precision, which allows for an econometric analysis at the highest resolution possible. 

\begin{figure}
\centering
\hspace*{-4ex}\hfill\subfigure[2:40:00 pm  -- 2:40:10 pm]{\includegraphics[height=0.35\textwidth]{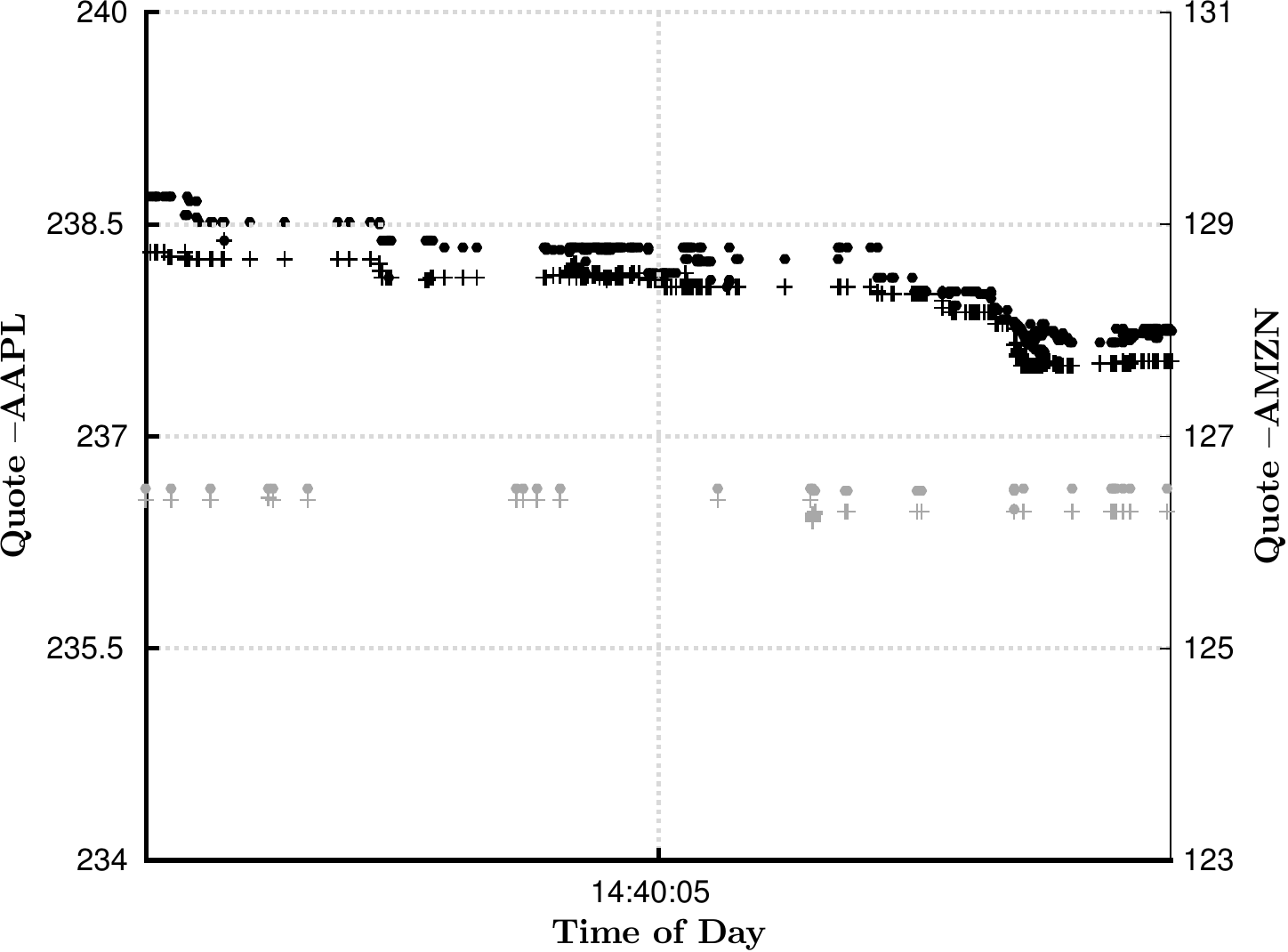}}\hspace*{-0.3ex}\hfill
\subfigure[2:50:00 pm  -- 2:50:10 pm]{\includegraphics[height=0.35\textwidth]{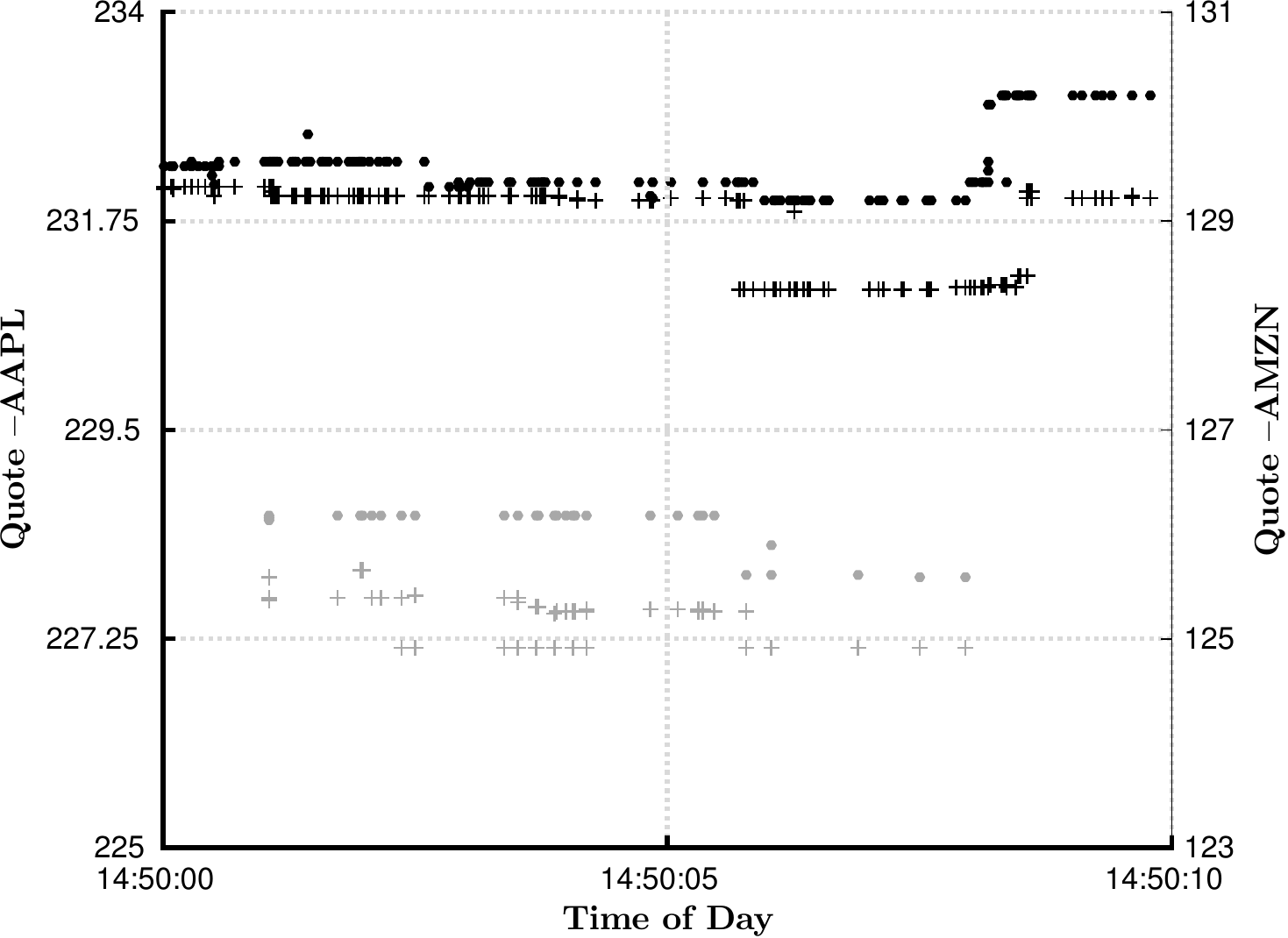}}\hfill
\caption{Bid- and ask-quotes of AAPL and AMZN on May 6th, 2010. Circles correspond to ask-quotes, crosses to bid-quotes. Black symbols correspond to AAPL, grey symbols to AMZN.}
\label{fig:quotes_AAPL_AMZN_20100605}
\end{figure}

In the web appendix, we provide summary statistics of the underlying raw data corresponding to the best ask and bid quotes in the Nasdaq market, recorded whenever the first level of the order book is updated. The average daily number of (level one) order book updates is $185,000$, corresponding to a new observation every $0.2$ seconds. Figure~\ref{fig:quotes_AAPL_AMZN_20100605} depicts the ask and bid quotes for Apple (AAPL) and Amazon (AMZN) during two ten-second time intervals at 2:40~pm and 2:50~pm on May 6, 2010. We observe three features: (i) most of the underlying order book updates do not cause a change in the best ask and bid quotes. Consequently, most of the corresponding event-to-event \emph{mid}-quote returns are zero. (ii) Despite the high precision of the time stamps, we observe time periods of several seconds without any activity, followed by other periods where order book activity is strongly clustered. Hence, the data is highly irregularly spaced. (iii) Orderbook updates of both stocks occur asynchronously over time. This is particularly due to obvious differences in event intensities. (iv) Particularly the Amazon quotes tend to bounce between different price levels. This is caused by considerable quoting activity and an obviously thin limit order book on the first level. The aforementioned effect even leads to a certain bouncing behavior in the resulting mid-quote returns, which is not necessarily attributed to movements of the underlying fundamental price, but rather to liquidity-induced noise.

\begin{figure}
\centering
\subfigure[Mid-quotes: 2:00 pm  -- 3:00 pm]{\includegraphics[height=0.35\textwidth]{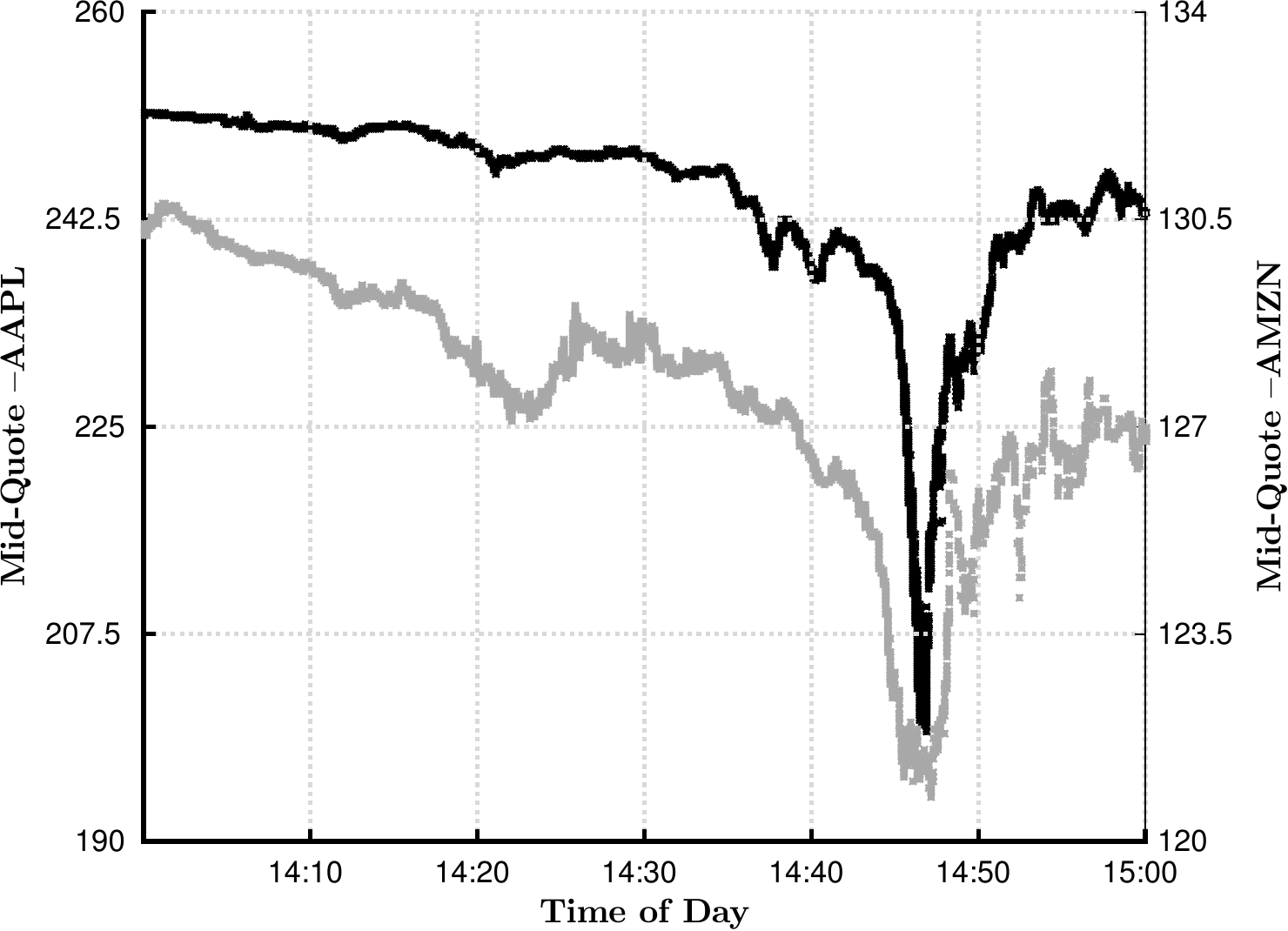}}\hfill
\subfigure[Spot correlations]{\includegraphics[height=0.35\textwidth]{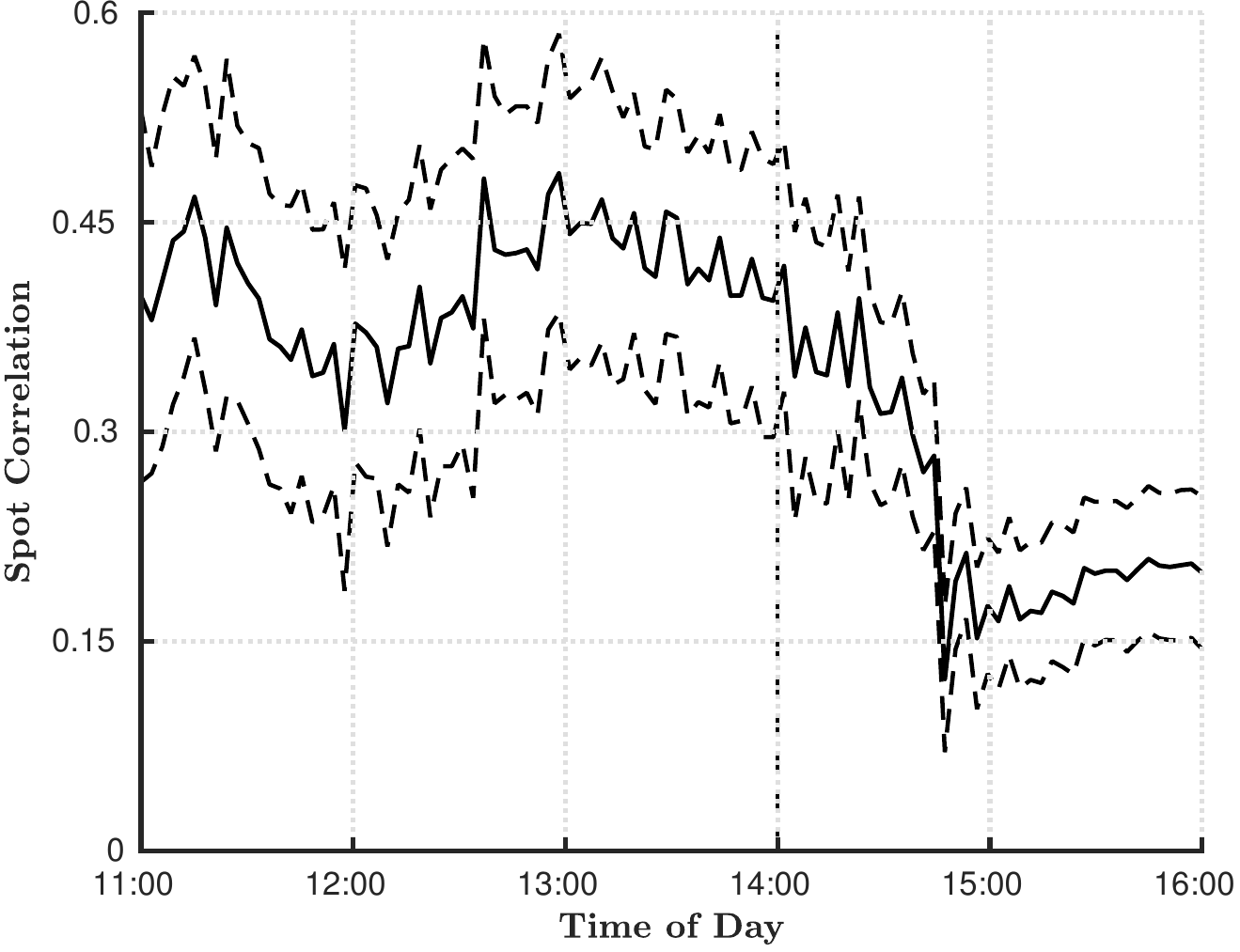}}\hfill
\caption{Mid-quotes and spot correlations of AAPL and AMZN on May 6th, 2010. In the left plot, black dots correspond to AAPL, grey dots to AMZN. In the right plot, dashed lines correspond to approximate pointwise $95\%$ confidence intervals according to Corollary~\ref{fclt} in Section \ref{sec:estimspotcov}.}
\label{fig:midq_corr_AAPL_AMZN_20100506}
\end{figure}

The left part of Figure~\ref{fig:midq_corr_AAPL_AMZN_20100506}, focusing on \emph{mid}-quotes, shows that the two above ten-second periods are at the heart of the flash crash occurring on May 6th, 2010 between 2:00~pm and 3:00~pm. While both figures reveal a microscopic and macroscopic view of the price behavior around this time point, neither picture provides hints on the underlying covariance and correlation between the two stocks and how they may change in such an extreme period.  The right part of Figure~\ref{fig:midq_corr_AAPL_AMZN_20100506} shows the behavior of the estimated correlation path during this trading day, revealing a strong and highly significant downward movement which results in a significantly different correlation level after the flash crash. The (approximate) confidence intervals depicted in this figure are constructed based on quote data for the two assets under focus on the given trading day only and rely on a feasible central limit theorem provided in this paper (see Section~\ref{sec:lmm_asymp}).

\begin{figure}
\centering
\hspace*{-4ex}\hfill\subfigure[1:11:00 pm  -- 1:11:10 pm]{\includegraphics[height=0.35\textwidth]{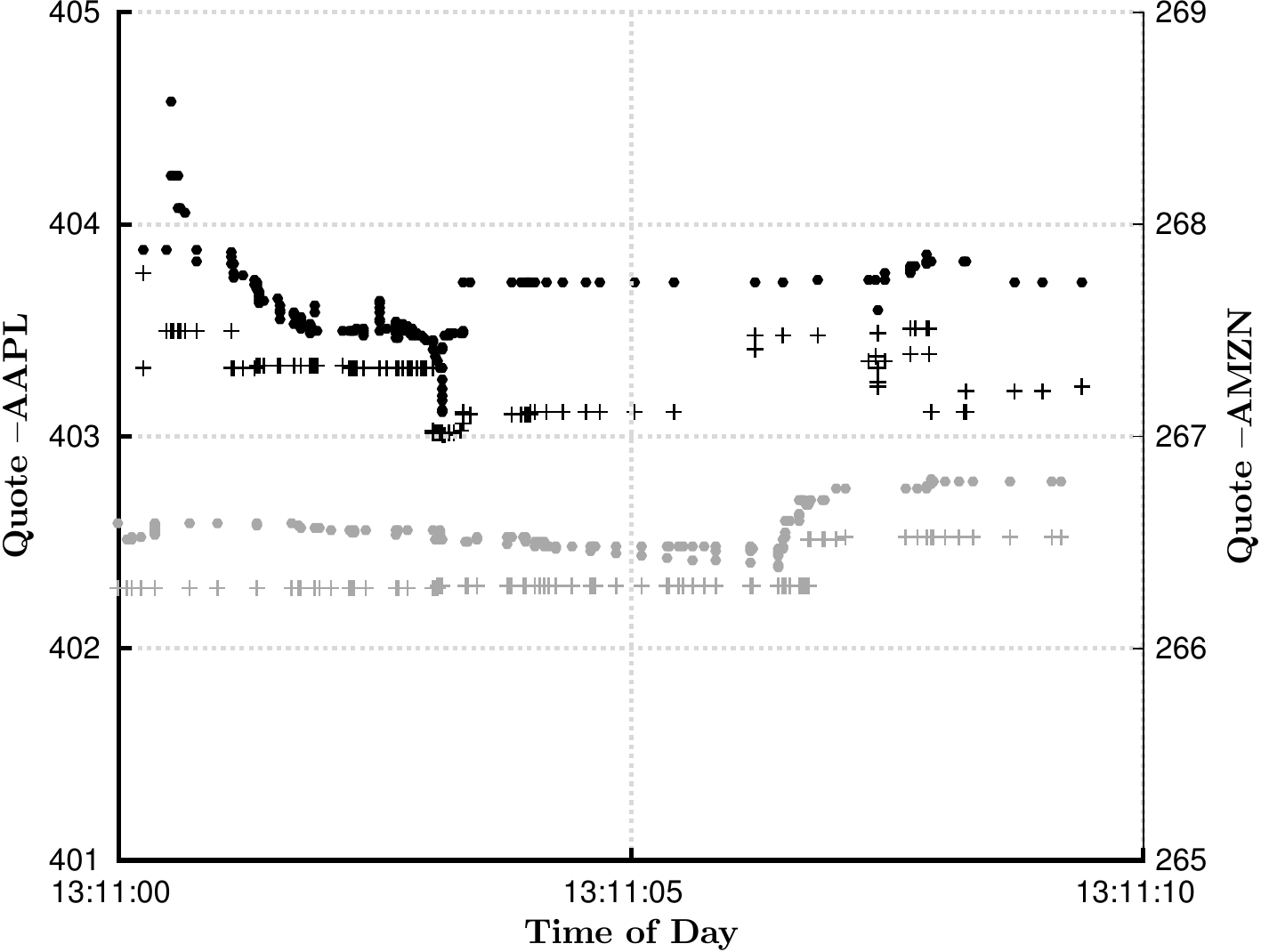}}\hspace*{-0.3ex}\hfill
\subfigure[1:13:00 pm  -- 1:13:10 pm]{\includegraphics[height=0.35\textwidth]{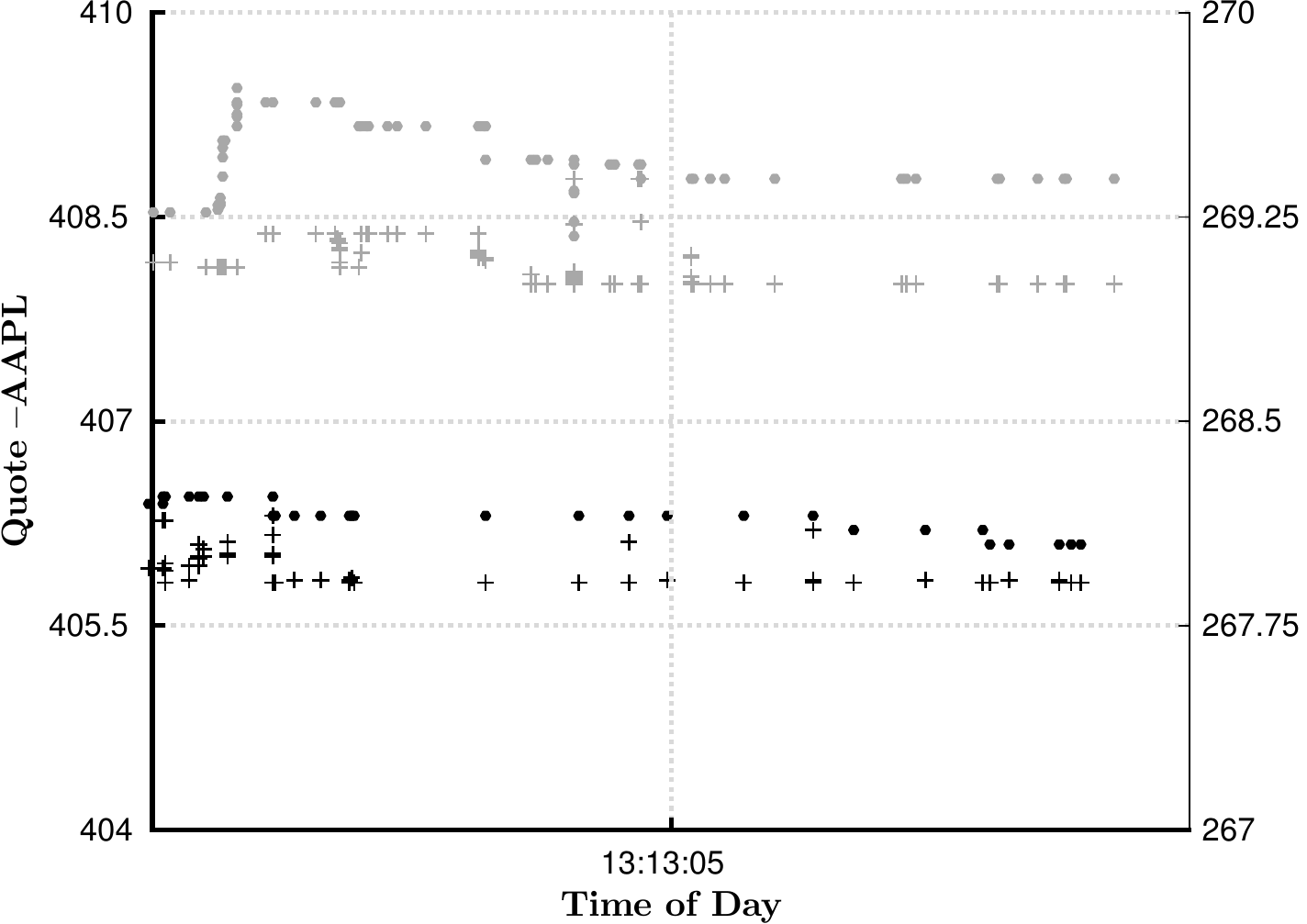}}\hfill
\caption{Bid- and ask-quotes of AAPL and AMZN on April 23rd, 2013. Circles correspond to ask-quotes, crosses to bid-quotes. Black symbols correspond to AAPL, grey symbols to AMZN.}
\label{fig:quotes_AAPL_AMZN_20130423}
\end{figure}

\begin{figure}
\centering
\subfigure[Mid-quotes: 1:05 pm  -- 1:15 pm]{\includegraphics[height=0.35\textwidth]{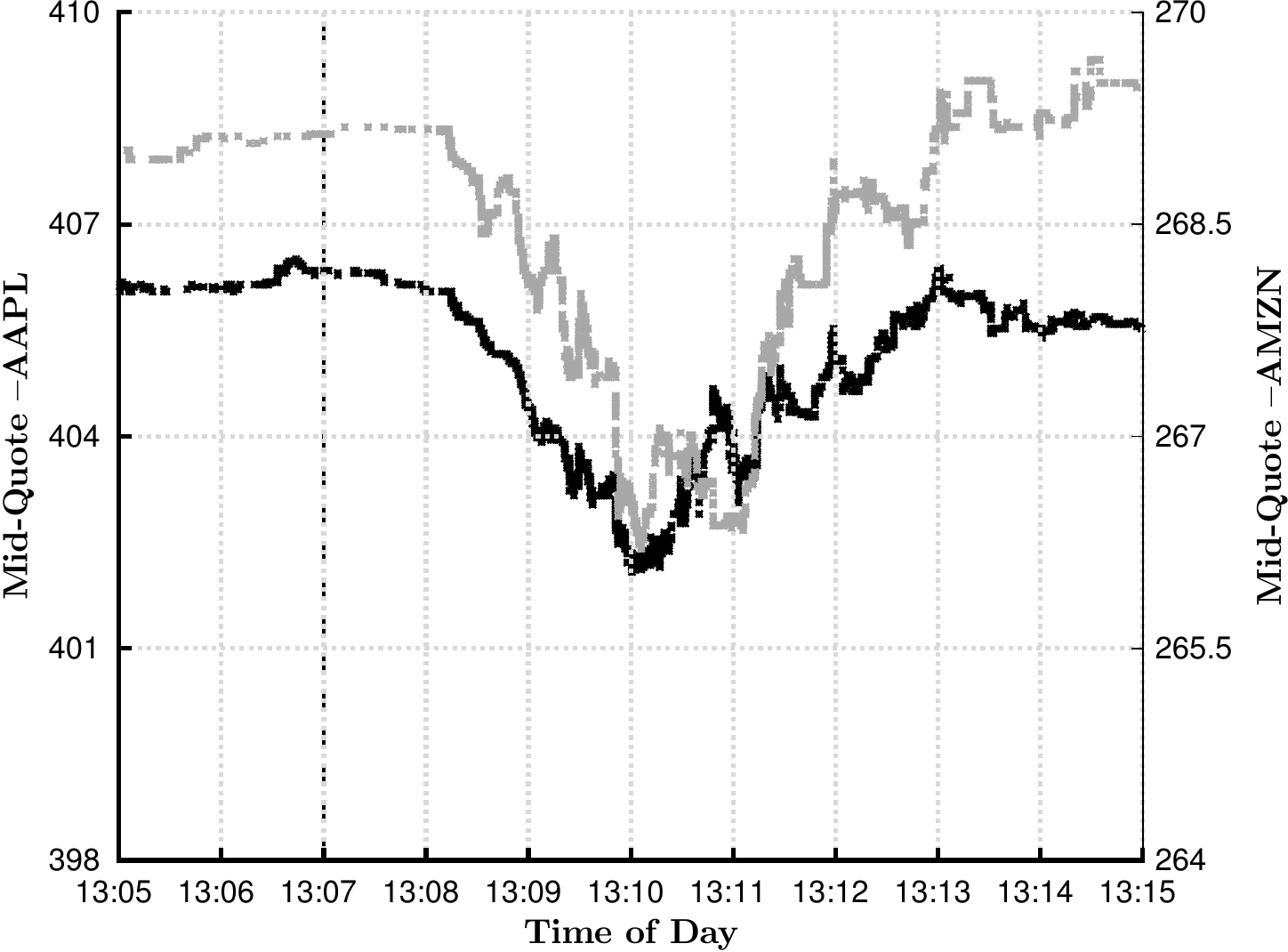}}\hfill
\subfigure[Spot correlations]{\includegraphics[height=0.35\textwidth]{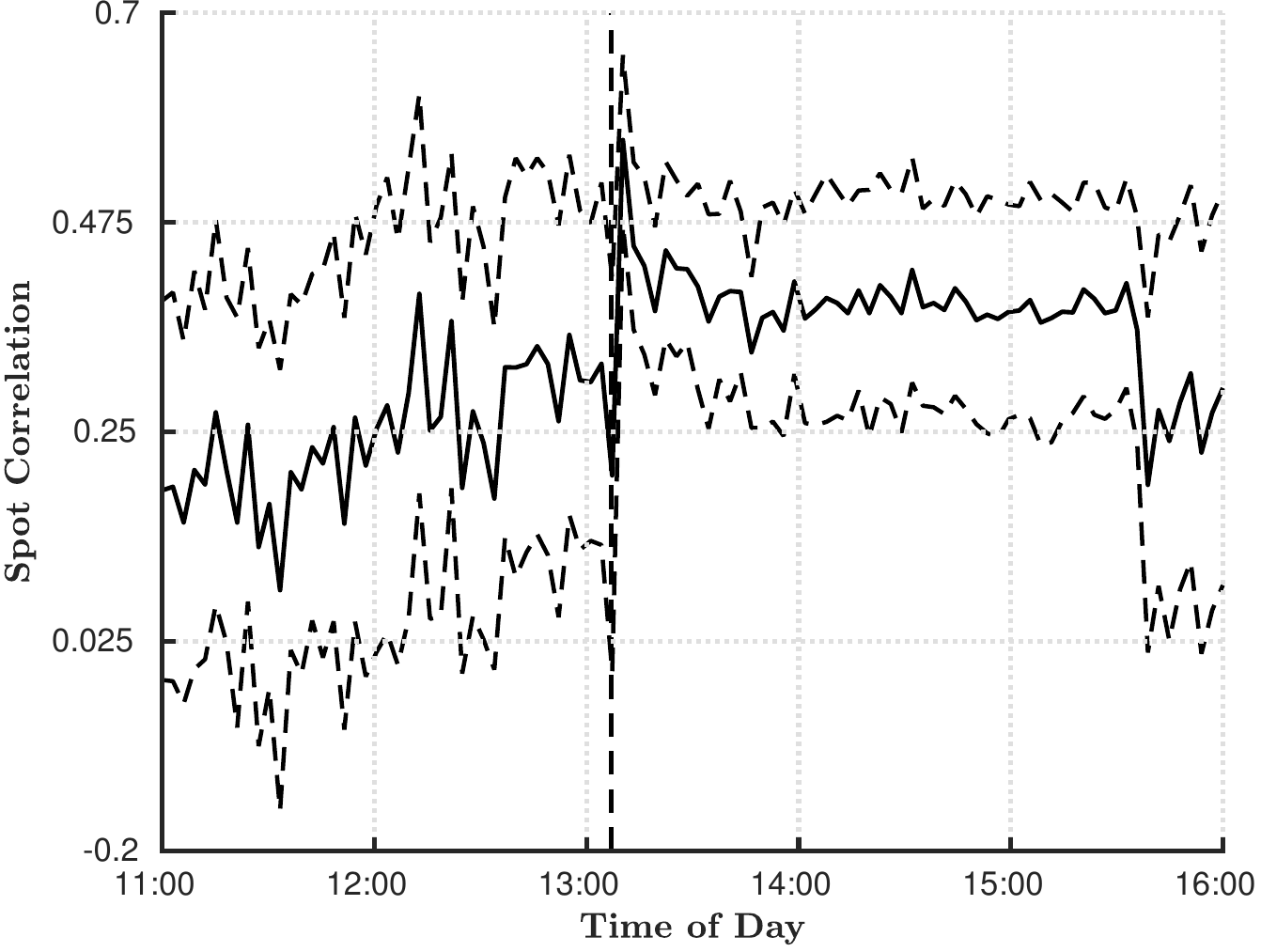}}\hfill
\caption{Mid-quotes and spot correlations of AAPL and AMZN on April 23rd, 2013. In the left plot, black dots correspond to AAPL, grey dots to AMZN. In the right plot, dashed lines correspond to approximate pointwise $95\%$ confidence intervals according to Corollary~\ref{fclt} in Section \ref{sec:estimspotcov}.}
\label{fig:midq_corr_AAPL_AMZN_20130423}
\end{figure}

Figure~\ref{fig:quotes_AAPL_AMZN_20130423} shows the ask/bid-quote behavior during two ten-second periods at 1:11~pm and 1:13~pm on April 23rd, 2013. While both pictures confirm the high-frequency properties of quote data discussed above, the asynchronicity of the two series becomes even more visible. As discussed in more detail in Section~\ref{sec:event} and illustrated in the left part of Figure~\ref{fig:midq_corr_AAPL_AMZN_20130423}, during this period, most equity prices dropped sharply because of faked Twitter news.  The right part of Figure~\ref{fig:midq_corr_AAPL_AMZN_20130423}  depicts the correlation path for this day, providing striking evidence for a strong and significant temporal shift in correlations. 


These two examples demonstrate that intraday movements in covariances and correlations can be substantial and can occur rapidly if new information arrives on the market. These movements can be empirically identified with sufficient precision, revealing important information for market surveillance and market microstructure research. The construction of correlation path estimates and corresponding confidence intervals, however, requires the optimal use of the underlying high-frequency information. From the illustrations above, it is obvious that a sufficiently precise identification of a correlation estimate at a single point in time (as, e.g., at 1:11~pm on April 23rd, 2013) cannot exploit information during the corresponding interval only, but needs to incorporate quote information from neighboring intervals. We therefore need to address the question of optimal smoothing over time, and thus the tradeoff between bias and variance. A further challenge is to account for the asynchronicity and irregular spacing of the observations, avoiding downward biases of estimates due to the Epps effect.
This is even more true as for some stocks a considerable amount of mid-quote returns equals zero and therefore does not necessarily provide new price information. In the web appendix, we show that on average only $13 \%$ of all mid-quote returns are non-zero. In the extreme case (e.g. for Microsoft), this quantity can amount to only  $1 \%$.  One initial step to utilize the underlying information in a (computationally) more efficient way is therefore to make use of quote \emph{revisions} only.  

Moreover, as sufficiently precise spot correlation estimates require exploiting high-frequency data on the highest possible frequency, correlation estimates need to be robust to possible market microstructure noise, i.e., deviations of the observed mid-quote price from the underlying ``true'' price process. Based on estimates of the (long-run) noise variance relying on an estimation procedure described in Section 1.3 of the web appendix and employing quote revisions, Table 5 of the the web appendix reports an average noise-to-signal ratio per observation of $1.5$. On ultra-high observation frequencies, market microstructure noise, however, is moreover likely to be serially correlated with the order of serial dependence being unknown ex-ante. In fact, using a test for serial correlation in the noise process as developed in the web appendix, the latter provides evidence for serial dependence up to an order of $12$ on average for mid-quote revisions.

\begin{figure}
\centering
\includegraphics[width=0.52\textwidth]{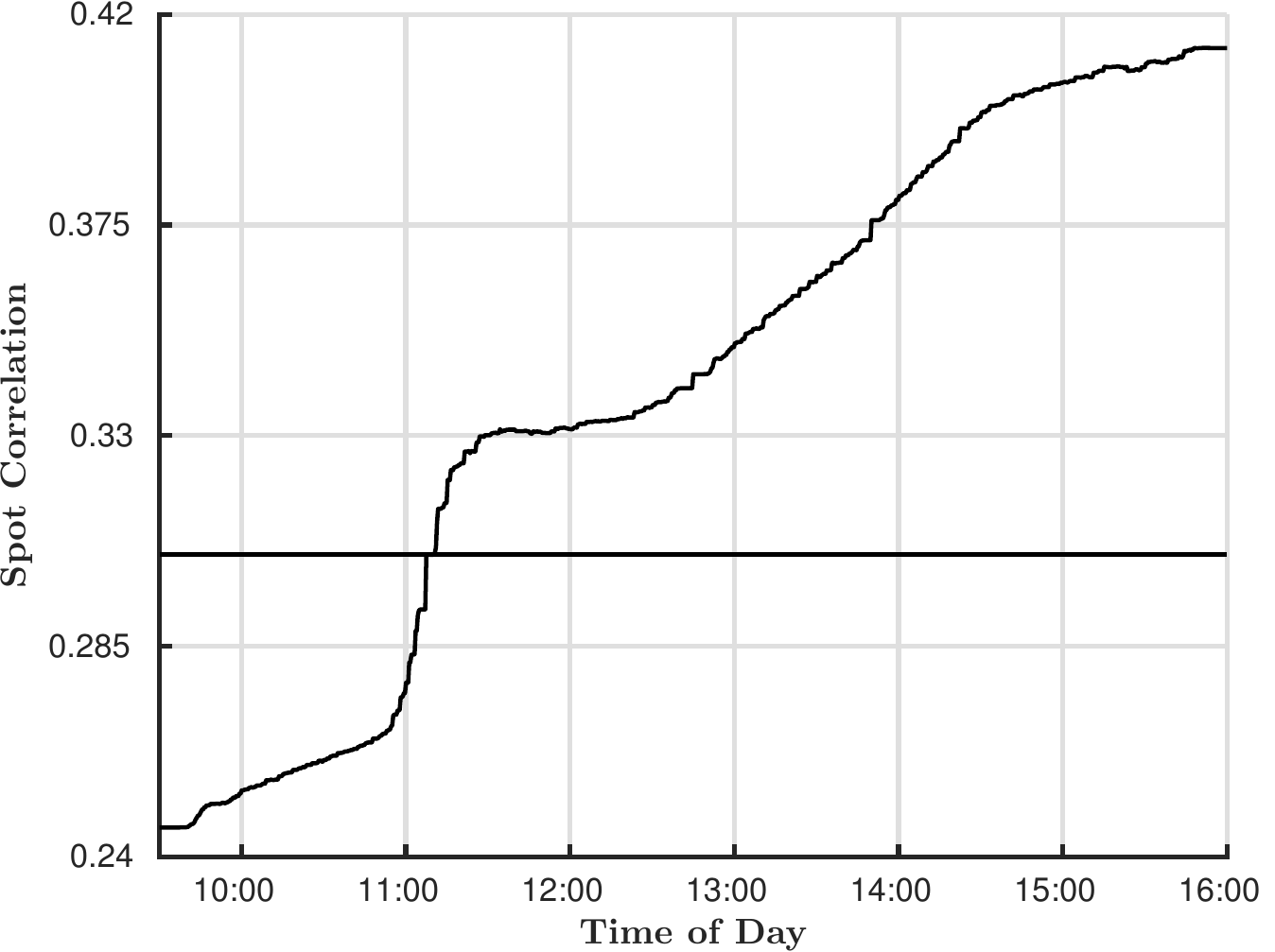}
\caption{Cross-sectional medians of across-day averages of spot correlations. Spot estimates are first averaged across days for each asset pair. Subsequently, cross-sectional sample medians of the across-day averages are computed. Solid horizontal line corresponds to the cross-sectional median of the across-day averages of \textit{integrated} correlation estimates. These are based on the LMM estimator of the integrated (open-to-close) covariance matrix by \citet{BHMR} accounting for serially dependent noise and using the same input parameter configuration as the spot estimators. 
``Unusual days'' discussed in Section~\ref{sec:event} as well as days with scheduled FOMC announcements are removed.}
\label{fig:med_corr_fulls}
\end{figure}

Finally, intraday trajectories of spot correlations are potentially subject to intraday periodicity effects. Indeed, one novel empirical finding of this paper is to identify distinct intraday seasonalities not only for individual asset return variances (as also documented in other work, e.g., in \citet{ande97,ande98}), but also for cross-asset correlations. Figure~\ref{fig:med_corr_fulls}  shows the cross-sectional medians of the across-day averages of pair-wise correlations, employing all combinations of the 30 most liquid Nasdaq stocks and all days through the period from May 2010 to April 2014 excluding ``unusual days'' as  discussed in Section~\ref{sec:event} as well as days with scheduled FOMC announcements. It turns out that correlations tend to systematically increase through the day with the highest rise during the morning hours.

To address these challenges and to incorporate these stylized facts of the data, it is thus necessary to construct an estimator which (i) optimally makes use of local mid-quote information as, e.g., depicted in Figures~\ref{fig:midq_corr_AAPL_AMZN_20100506} and~\ref{fig:midq_corr_AAPL_AMZN_20130423}, resulting in consistent and precise estimates with the highest convergence rates possible, (ii) allows for fairly general properties of the underlying spot volatility (matrix) process, incorporating, e.g., intraday periodicities, (iii) accounts for serially correlated and potentially endogenous noise, and (iv) yields feasible asymptotic inference, accounting for the pre-estimation of noise-dependent quantities. Based on the asymptotic results for such a spot covariance matrix estimator and the Delta-method, consistent estimators for other quantities of interest, such as spot betas and spot correlations can be deduced.

\section{Estimation of Spot Covariances\label{sec:estimspotcov}}
\subsection{Theoretical Setup and Assumptions\label{sec:setup}}

Let $\left(X_t\right)_{t\ge 0}$ denote the $d$-dimensional efficient log-price process. In line with the literature and motivated by well-known no-arbitrage arguments, we assume that $X_t$ follows a continuous It\^{o} semi-martingale
\begin{align}\label{sm}
X_t=X_0+\int_0^t b_s\,ds+\int_0^t\sigma_s\,dB_s,\;t\in[0,1],
\end{align}
defined on a filtered probability space $\bigl(\Omega,\mathcal{F},\left(\mathcal{F}\right)_{t\ge 0},\P\bigr)$ with drift $b_s$, $d$-dimensional standard Brownian motion $B_s$ and instantaneous volatility matrix $\sigma_s$. The latter yields the $(d\times d)$-dimensional spot covariance matrix $\Sigma_s=\sigma_s\sigma_s^{\top}$, which is our object of interest. We consider a setting in which discrete and non-synchronous observations of the process~$\eqref{sm}$ are diluted by market microstructure noise, i.e.,
\begin{align}
\label{obsp}
Y_i^{(p)}=X_{t_i^{(p)}}^{(p)}+\epsilon_i^{(p)}\,,\,i=0,\ldots,n_p,\;p=1,\ldots,d\,,
\end{align}
with observation times $t_i^{(p)}$, and observation errors $\epsilon_i^{(p)}$. Observed returns for component $p\in\{1,\ldots,d\}$ are given by
\begin{align}
\label{obs}
\Delta_i Y^{(p)}&=Y_i^{(p)}-Y_{i-1}^{(p)}=\Delta_i X^{(p)}+\Delta_i\epsilon^{(p)}\\\nonumber
&=X_{t_i^{(p)}}^{(p)}-X_{t_{i-1}^{(p)}}^{(p)}+\epsilon_i^{(p)}-\epsilon_{i-1}^{(p)},\quad i=1,\ldots,n_p.
\end{align}

Let $n=\min\limits_{1\leq p\leq d} n_p$ denote the number of observations of the least liquid asset. In Section~\ref{sec:lmm_asymp}, we consider high-frequency asymptotics with $n/n_p\rightarrow \nu_p$ for constants $0<\nu_p\leq 1$, such that the asymptotic variance-covariance matrices for estimators of $\Sigma_s$ are regular. 

Below we summarize the assumptions on the instantaneous volatility matrix and drift process, noise properties and observation times. In order to describe smoothness classes of the spot covariance matrix process and other functions, we consider balls in H\"older spaces of order $\alpha\in(0,1]$ and with radius $R>0$:
\begin{align*}
&C^{\alpha,R}([0,1],\mathcal{E})=\{f:[0,1]\rightarrow \mathcal{E}|\;\|f\|_{\alpha}\le R\}\,,\,\|f\|_{\alpha}\:=\|f\|_{\infty}+\sup_{x\ne y}{\frac{\|f(x)-f(y)\|}{|x-y|^{\alpha}}}\,,
\end{align*}
where  $\|\cdot\|$ denotes the usual spectral norm and 
$\|f\|_{\infty}\:=\sup_{t\in[0,1]}{\|f(t)\|}$ for functions on $[0,1]$. In our setup, we have $\mathcal{E}=\mathds{R}^{d\times d'}$ for matrix-valued functions, $\mathcal{E}=\mathds{R}^d$ for vectors or $\mathcal{E}=[0,1]$ for distribution functions.

First, for the drift process in~\eqref{sm}, we only assume a very mild regularity:
\begin{sassump}
\label{drift}
$(b_s)_{s\in[0,1]}$ is an ($\mathcal{F}_s$)-adapted process with $b_s=g(b_s^{(1)},b_s^{(2)})$, $g:\R^{2d}\rightarrow\R^{d}$ being a continuously differentiable function in all coordinates, $b_s^{(1)}$ an It\^{o} semi-martingale with locally bounded characteristics and $b_s^{(2)}\in C^{\nu,R}([0,1],\mathds{R}^d)$ for some $R<\infty$ and some $ \nu>0$.
\end{sassump}
Assumptions on the instantaneous volatility matrix process in~\eqref{sm} can be summarized as:
\begin{sassump}
\label{sigma}
(i) $(\sigma_s)_{s\in[0,1]}$ follows an ($\mathcal{F}_s$)-adapted process satisfying $ \Sigma_s=\sigma_s\sigma_s^{\top}\ge \underline{\Sigma}$ uniformly for some strictly positive definite matrix $\underline{\Sigma}$.

 (ii)  $(\sigma_s)_{s\in[0,1]}$ satisfies $\sigma_s=f\big(\sigma_s^{(1)},\sigma_s^{(2)}\big)$ with $f:\R^{2d\times 2d'}\rightarrow\R^{d\times d'}$ being a continuously differentiable function in all coordinates, where
\begin{itemize}
\item For $\alpha\in (0,1/2]$, $\sigma_s^{(1)}$ is an It\^{o} semi-martingale with locally bounded characteristics.
\item
$\sigma_s^{(2)}\in C^{\alpha,R}\big([0,1],\mathds{R}^{d\times d'}\big)$ with some $R<\infty$.
\end{itemize}
\end{sassump}
Hence, $\sigma_s$ is a function of an It\^{o} semi-martingale $\sigma_s^{(1)}$ and an additional H\"older smooth component $\sigma_s^{(2)}$. The latter can capture intraday periodicity effects \citep[see, e.g.,][]{ande97}.
Assumption \ref{sigma} depends on the smoothness parameter $\alpha$ and reads similar as Assumption (K-$v$) in \cite{jacodtodorov}. The larger $\alpha$, the more restrictive becomes Assumption \ref{sigma}. If $\alpha > 1/2$, we assume that the semi-martingale component $\sigma_s^{(1)}$ vanishes and $\sigma_s$ is exclusively driven by the component $\sigma_s^{(2)}$. Hence, the more interesting case is $\alpha\le 1/2$. Then, Assumption \ref{sigma} allows also for a semi-martingale volatility with volatility jumps. Importantly, the above assumptions also allow for leverage effects, that is, a non-zero correlation between $\sigma_s$ and the Brownian motion $B_s$ in~\eqref{sm}. It is natural to develop results under this general smoothness assumption depending on $\alpha$ as it is commonly known that in nonparametric estimation problems, the underlying regularity $\alpha$ determines the size of smoothing windows and a fortiori the resulting (optimal) convergence rates.

Our assumptions on the microstructure noise process in~\eqref{obsp} are stated in observation time, which is in line with, e.g., \citet{hans06} and \citet{bn1}:
\begin{sassump}
\label{noise}
(i) $\epsilon=\{\epsilon_i^{(p)},i=0,\ldots,n_p,p=1,\ldots,d\}$ is independent of $X$ and has independent components, i.e.,\,$\epsilon_i^{(p)}$ is independent of $\epsilon_j^{(q)}$ for all $i,j$ and $p\ne q$.

(ii) At least the first eight moments of $\epsilon_i^{(p)},i=0,\ldots,n_p,$ exist for each $p=1,\ldots,d$.

(iii) $\epsilon_i^{(p)},i=0,\ldots,n_p,$ follows an $R$-dependent process for some $R<\infty$, implying that $\cov\big(\epsilon_i^{(p)},\epsilon_{i+u}^{(p)}\big)=0$ for $u>R$ and each $p=1,\ldots,d$. Define by
\begin{align}\label{eta}\eta_p=\eta_0^{(p)}+2\sum_{u=1}^R\eta_u^{(p)},~~\mbox{with}~~\eta_u^{(p)}\:=\cov\big(\epsilon_i^{(p)},\epsilon_{i+u}^{(p)}\big),u\le R,\end{align}
the component-wise long-run noise variances, where the $\eta_u^{(p)},0\le u\le R$, are constant for all $0\le i\le n-u$. We impose that $\eta_p>0$ for all $p$.
\end{sassump}
The independence between noise and the efficient price, as stated in part~(i) of Assumption~\ref{noise}, is standard in the literature \citep[see, e.g.,][]{zhangmykland}. On the other hand, however, \cite{hans06} report evidence for an endogeneity with dependence between noise and efficient price and it is of interest that estimators are robust in that case. To meet this objective, we show in Section 1.3 of the web appendix that our estimator will be robust, in the sense of keeping the same asymptotic properties, to correlations between signal and noise. Considering serially dependent noise is non-standard and motivated by empirical results, e.g., in \citet{hans06}. The moving-average-type dependence structure in the noise process in part~(iii) of Assumption~\ref{noise} follows, e.g., \citet{haut13}, implying the long-run variance \eqref{eta}. Generalized realized kernel estimation of the covariation in a very general model with endogenous and serially correlated noise has been presented in \cite{varneskov}. Cross-sectional dependence of the noise is left aside in our theory, since a notion of simultaneous dependence in the presence of non-synchronicity is by now not established. In principle, non-diagonal noise variance-covariance matrices can be included in the multivariate framework below, see \cite{spectral} for a setup allowing for cross-sectional dependence when recordings are synchronous.

Finally, we assume that the timing of observations in~\eqref{obsp} is driven by c.d.f.'s~$F_p$ governing the transformations of observation times to equidistant sampling schemes by means of suitable quantile transformations:
\begin{sassump}
\label{obs_ass}
There exist differentiable cumulative distribution functions $F_p, p=1,\ldots,d$, such that the observation regimes satisfy $t_i^{(p)}=F_p^{-1}(i/n_p)$ , $0\le i\le n_p,p\in\{1,\ldots,d\}$, where $F_p'\in C^{\alpha,R}\big([0,1],$ $[0,1]\big),p=1,\ldots,d$, with $\alpha$ being the smoothness exponent in Assumption \ref{sigma} for some $R<\infty$. $F_p,p=1,\ldots,d$, can be random, but independent from the observed process $\{Y_i^{(p)}\}$.
\end{sassump}
A treatment of endogenous times in the given theoretical framework is beyond the scope of this paper. See \citet{koike} for a recent study of endogenous times and \citet{limyk14} for a study in a setting neglecting microstructure noise.\\
Combining time-invariant (long-run) noise variances $\eta_p$ and locally different observation frequencies from Assumptions~\ref{noise} and~\ref{obs_ass} implies locally varying noise levels. In the asymptotic framework with $n/n_p\rightarrow \nu_p$, where $0<\nu_p\le 1,p=1,\ldots,d$, for $n\rightarrow\infty$, we define the continuous-time noise level matrix
\begin{align}\label{noiselevel}{H}_s=\operatorname{diag}\big((\eta_p\nu_p(F_p^{-1})^{\prime}(s))^{1/2}\big)_{1\le p\le d}\,.\end{align}
Note that for equally-spaced observations, we have $F_p(s)=s$, such that $(F_p^{-1})^{\prime}(s)=1$. Then, the $p$-specific (asymptotic) noise level is $(\eta_p\nu_p)^{1/2}$ with the constant $\nu_p$ expressing the inverse of the sample size of the $p$-th process relative to the ``least liquid'' process. Hence,  having less frequent observations on a sub-interval is equivalent to having higher noise dilution by microstructure effects on this sub-interval. This interplay between noise and liquidity has been discussed by \citet{BHMR}.

\subsection{Local Method of Moments Estimation of the Spot  Covariance Matrix\label{sec:lmm}}

Our approach for estimating the instantaneous covariance matrix rests upon the concept of the local method of moments (LMM) introduced in \cite{BHMR}. We partition the interval $[0,1]$ into equidistant blocks $[kh_n,(k+1)h_n],k=0,\ldots,h_n^{-1}-1$, with the block length $h_n$ asymptotically shrinking to zero, $h_n\rightarrow 0$ as $n\rightarrow\infty$. The key idea is to approximate the underlying process \eqref{sm} in model \eqref{obsp} by a process with block-wise constant covariance matrices and noise levels. In the (more simplified) setting of \cite{BHMR}, it is shown that such a locally constant approximation induces an estimation error for the integrated covariation, which, however, can be asymptotically neglected for sufficient smoothness of $\Sigma_t$ and $F_p$ if the block sizes $h_n$ shrink sufficiently fast with increasing $n$. This opens the path to construct an asymptotically efficient estimator of the integrated covariation matrix based on optimal block-wise estimates.

In the present setting, we build on the idea of block-wise constant approximations of the underlying covariance and the noise process and show that it allows constructing a consistent spot covariance estimator, which can attain an optimal rate. A major building block is the construction of an unbiased estimator of the block-wise covariance matrix $\Sigma_{k h_n}=\sigma_{k h_n}\sigma_{k h_n}^{\top}$ based on the local spectral statistics
\begin{align}\label{spec}
S_{jk}=\pi jh_n^{-1}\Bigg(\sum_{i=1}^{n_p}\Big( Y_{i}^{(p)}-Y_{i-1}^{(p)}\Big)\Phi_{jk}
\Big(\frac{t_{i-1}^{(p)}+t_{i}^{(p)}}{2}\Big)\Bigg)_{1\le p\le d}\,,
\end{align}
where $\Phi_{jk}$ denote orthogonal sine functions with (spectral) frequency $j$, whose derivatives $\Phi_{jk}'$ form another orthogonal system corresponding to the eigenfunctions of the covariance operator of a Brownian motion, and are given by
\begin{align}\label{Phi}
\Phi_{jk}(t)&=\frac{\sqrt{2h_n}}{j\pi} \sin{\left(j\pi h_n^{-1}\left(t-kh_n\right)\right)}\1_{[kh_n,(k+1)h_n)}(t),j\ge 1\,.
\end{align}
The statistics~\eqref{spec} de-correlate the noisy observations~\eqref{obs} and can be thought of as representing their block-wise principal components. They bear some resemblance to the pre-averaged returns as employed in \citet{JLMPV}. While pre-averaging estimators, however, utilize rolling (local) windows around each observation, our approach relies on fixed blocks and optimal combinations in the spectral frequency domain. Related approaches for a univariate framework can be found in \citet{hans08}, \citet{reiss} and \citet{corsi}. 

It has been shown in \cite{stable} that
\begin{align} \label{eq_distr_S}
\cov(S_{jk}) = (\Sigma_{kh_n}+\pi^2j^2h_n^{-2} {\bf{H}}_k^n)(1+ \KLEINO(1)),
\end{align}
where ${\bf H}_k^n$ denotes the block-wise constant diagonal noise level matrix with entries
\begin{align}\label{noiselevel2}
\big({\bf{H}}_k^n\big)^{(pp)}=n_p^{-1}\eta_p (F_p^{-1})'(kh_n)\,.
\end{align}
The relation \eqref{eq_distr_S} corresponds to Equation (2.4) in \cite{BHMR} plus a negligible remainder in the more general model and suggests estimating $\Sigma_{kh_n}$ based on the empirical covariance $S_{jk}S_{jk}^{\top}$, which is bias-corrected by the noise-induced term $\pi^2j^2h_n^{-2} {\bf{H}}_k^n$.

An initial (pre-) estimator of the spot covariance matrix at time $s\in[0,1]$, $\Sigma_{s}$, is then constructed based on  bias-corrected block-wise empirical covariances $S_{jk}S_{jk}^{\top}$, which are averaged across spectral frequencies $j=1,\ldots, J_{n}^{p}$, and  a set of adjacent blocks,
\begin{align}
\label{pilot}
\operatorname{vec}\big(\hat\Sigma_{kh_n}^{pre}\big)=(U_{s,n}-L_{s,n}+1)^{-1}\sum_{k=L_{s,n}}^{U_{s,n}}{\left(J_{n}^{p}\right)}^{-1}\sum_{j=1}^{J_{n}^{p}}\operatorname{vec}\Big(S_{jk}S_{jk}^{\top}-\pi^2j^2h_n^{-2}\hat{\bf{H}}_k^n\Big)\,,
\end{align}
with  $L_{s,n}=\max\{\lfloor sh_n^{-1}\rfloor -K_{n}, 0\}$ and $U_{s,n}=\min\{\lfloor sh_n^{-1}\rfloor +K_{n}, \lceil h_n^{-1}\rceil-1\}$  for a two-sided estimator as well as $L_{s,n}=\max\{\lfloor sh_n^{-1}\rfloor -2 K_{n}, 0\}$ and $U_{s,n}=\min\{\lfloor sh_n^{-1}\rfloor, \lceil h_n^{-1}\rceil-1\}$ for a one-sided estimator, such that the length of the smoothing window obeys $U_{s,n}-L_{s,n}+1 \le 2K_{n}+1$. In this context, two-sided means that at some time $s$, we estimate the covariances by locally smoothing over a window centered around $s$. One-sided refers to the same method, but with smoothing over a window before and up to time $s$. Asymptotic properties are the same.
$\hat{\bf{H}}_k^n$ is a $\sqrt{n}$-consistent estimator of ${\bf{H}}_k^n$ with $p$-th diagonal element
\begin{align}
\label{noiselevelest}
\big(\hat{\mathbf{H}}_k^n\big)^{(pp)}=\frac{\hat\eta_p}{h_n}\sum_{kh_n\le t_i^{(p)}\le (k+1)h_n}\big(t_i^{(p)}-t_{i-1}^{(p)}\big)^2.
\end{align}
Details on the construction of the estimator of the component-wise long-run noise variances, $\hat\eta_p$, are provided in Section 1.3 of the web appendix.

For each spectral frequency $j$, the statistic $S_{jk}S_{jk}^{\top}-\pi^2j^2h_n^{-2}\hat{\bf{H}}_k^n$ is an (asymptotically) unbiased though inefficient estimator of $\Sigma_{k h_n}$. Averaging across different frequencies therefore increases the estimator's efficiency. Equally weighting as in \eqref{pilot}, however, is not necessarily optimal. A more efficient estimator can be devised by considering~\eqref{pilot} as the pre-estimated spot covariance matrix and then, derive estimated optimal weight matrices $\hat W_j$, yielding the final LMM spot covariance matrix estimator as
\begin{align}
\label{spotcov}
\operatorname{vec}\big(\hat\Sigma_s\big)=(U_{s,n}-L_{s,n}+1)^{-1}\sum_{k=L_{s,n}}^{U_{s,n}}\sum_{j=1}^{J_{n}}\,&\hat{W}_{j}\big(\hat{\bf H}_k^n,\hat\Sigma_{kh_n}^{pre}\big)\\\nonumber
&\times\operatorname{vec}\Big(S_{jk}S_{jk}^{\top}-\pi^2j^2h_n^{-2}\hat{\bf{H}}_k^n\Big)\,.
\end{align}
As outlined in detail in \cite{BHMR}, the true optimal weights are given proportionally to the local Fisher information matrices according to
\begin{align}
\label{weights}
W_{j}\big({\bf{H}}_k^n,\Sigma_{kh_n}\big)&=\Big(\sum_{u=1}^{J_n}\big(\Sigma_{kh_n}+\pi^2u^2h_n^{-2}{\bf{H}}_k^n\big)^{-\otimes 2}\Big)^{-1}\big(\Sigma_{kh_n}+\pi^2j^2h_n^{-2}{\bf{H}}_k^n\big)^{-\otimes 2} \\
&=I_k^{-1}I_{jk}, \nonumber
\end{align}
with $I_{jk}$ being the Fisher information matrix associated with block $k$ and spectral frequency $j$, given by
\begin{align}
\label{fisher_info}
I_{jk}=\tfrac12 \big(\Sigma_{kh_n}+\pi^2j^2h_n^{-2}{\bf{H}}_k^n\big)^{-\otimes 2},
\end{align}
and $I_k=\sum_{j=1}^{J_n}I_{jk}$ denoting the $k$-specific Fisher information (exploiting the independence across frequencies $j$). Here, $A^{\otimes 2}=A\otimes A$ denotes the Kronecker product of a matrix with itself and $A^{-\otimes 2}=A^{-1}\otimes A^{-1}=(A\otimes A)^{-1}$. We show in Section~\ref{sec:lmm_asymp} that the estimator~\eqref{spotcov}, which builds on the idealized model considered in \citet{BHMR}, is consistent and satisfies a stable CLT under the more realistic and general assumptions of Section~\ref{sec:setup}.

While both the pilot estimator~\eqref{pilot} and the LMM estimator~\eqref{spotcov} are symmetric, neither is guaranteed to yield positive semi-definite estimates. Confidence is based on estimated Fisher information matrices $\hat I_k$, see \eqref{weights}, which are by construction positive-definite. For the estimates themselves, we can set negative eigenvalues equal to zero, which is tantamount to a projection on the space of positive semi-definite matrices. This adjustment does not affect the asymptotic properties of the estimator. For a similar adjustment to a realized kernel integrated covariance estimator see \cite{varneskov}.

\subsection{Asymptotic Properties\label{sec:lmm_asymp}}

As a prerequisite for the discussion of the central limit theorem for the estimator~\eqref{spotcov}, some considerations regarding $K_n$, which determines the length of the smoothing window, are needed. For this purpose, suppose that a certain smoothness $\alpha\in(0,1]$ of the instantaneous volatility matrix is granted according to Assumption~\ref{sigma}. Then, a simple computation yields $\|\Cov\big(\hat\Sigma_s\big)\|=\mathcal{O}\big(K_n^{-1}\big)$, implying a bias-variance trade-off in the mean square error $\operatorname{MSE}\big(\hat\Sigma_s\big):=\E\left[\|\hat\Sigma_s-\Sigma_s\|^2\right]$. More precisely, for a specific $\alpha>0$, we have 
\begin{align}\label{mse}\operatorname{MSE}\big(\hat\Sigma_s\big)=\mathcal{O}\big(K_n^{-1}\big)+\mathcal{O}\big(K_n^{2\alpha}h_n^{2\alpha}\big)\,,\end{align}
where the first term originates from the variance and the second term is induced by the squared bias. Consequently, for given $h_n\propto \log(n)n^{-1/2}$, which optimally balances noise and discretization error as derived in \cite{BHMR}, choosing $K_n\propto n^{\alpha/(2\alpha+1)}$ minimizes the MSE and facilitates an estimator with $\sqrt{K_n}$ convergence rate. Finally, the desired central limit theorem for the estimator~\eqref{spotcov} requires a slight undersmoothing, resulting in a smaller choice of $K_n$:
\begin{theo}
\label{cltspot}
We assume a setup with observations of the type \eqref{obsp}, a signal \eqref{sm} and the validity of Assumptions~\ref{drift}-\ref{obs_ass}. Then, for $h_n=\kappa_1 \log{(n)}n^{-1/2}$, $K_n=\kappa_2 n^{\beta}(\log{(n)})^{-1}$ with constants $\kappa_1,\kappa_2$ and $0<\beta<\alpha(2\alpha+1)^{-1}$, for $J_n\rightarrow\infty$ and $n/n_p\rightarrow\nu_p$ with $0<\nu_p\le 1,p=1,\ldots,d$, as $n\rightarrow\infty$, the spot covariance matrix estimator~\eqref{spotcov} satisfies the pointwise $\mathcal{F}$-stable central limit theorem:
\begin{align}
\label{cltspoteq}
n^{\beta/2}\operatorname{vec}\big(\hat \Sigma_s-\Sigma_s\big)\stackrel{d-(st)}{\longrightarrow}\mathbf{N}\Big(0,2\big(\Sigma\otimes \Sigma_{H}^{1/2}+\Sigma_{H}^{1/2}\otimes \Sigma\big)_s\,\mathcal{Z}\Big),\quad s\in[0,1]\,,
\end{align}
where $\Sigma_H=H\big(H^{-1}\Sigma H^{-1}\big)^{1/2}H$, with noise level $H$ from \eqref{noiselevel} and $\mathcal{Z}=\Cov(\operatorname{vec}(ZZ^{\top}))$ for $Z\sim \mathbf{N}(0,E_d)$ being a standard normally distributed random vector.
\end{theo}
Theorem~\ref{cltspot} is proved in the web appendix. Though Assumption \ref{sigma} involves volatility jumps for $\alpha\le 1/2$, \eqref{cltspoteq} applies, because for any fixed $s\in[0,1]$, the probability of a jump in the asymptotically small smoothing window converges to zero. For finite-sample applications of estimating the spot covariance matrix in the vicinity of a structural change, however, one should carefully adjust the chosen smoothing windows. Lemma~1 in the web appendix provides the key step to extend the analysis to autocorrelated noise. Its proof reveals, at the same time, why the generalization from Gaussian i.i.d.\,noise in \cite{BHMR} to the general Assumption \ref{noise} does not affect the asymptotic variance of the estimator. The convergence in~\eqref{cltspoteq} is $\mathcal{F}$-stable, which is equivalent to joint weak convergence with any $\mathcal{F}$-measurable bounded random variable defined on the same probability space as $X$. This allows for a feasible version of the limit theorem, even for general stochastic volatilities with leverage effects if we re-scale the estimator by the inclusively obtained estimated variance:
\begin{corro}
\label{fclt}
Under the assumptions of Theorem~\ref{cltspot}, the spot covariance matrix estimator~\eqref{spotcov} satisfies the feasible central limit theorem given by
\begin{subequations}
\begin{align}
\label{fclteq}
(U_{s,n}-L_{s,n}+1)^{1/2}\big(\hat{ \mathds{V}}_s^n\big)^{-1/2}\operatorname{vec}\big(\hat \Sigma_s-\Sigma_s\big)\stackrel{d}{\rightarrow}\mathbf{N}\Big(0,\mathcal{Z}\Big),\quad s\in[0,1]\,,
\end{align}
\begin{align}
\label{variance}
\text{where}\qquad\hat{ \mathds{V}}_s^n=(U_{s,n}-L_{s,n}+1)^{-1}\sum_{k=L_{s,n}}^{U_{s,n}}\Bigg(\sum_{j=1}^{J_n}\hat{I}_{jk}\Bigg)^{-1},\qquad\qquad\qquad
\end{align}
\end{subequations}
with $U_{s,n}$ and $L_{s,n}$ defined as in~\eqref{pilot} and \eqref{spotcov}. $\hat{I}_{jk}$ is defined according to~\eqref{fisher_info} with ${\bf{H}}_k^n$ and $\Sigma_{kh_n}$, $k=0,\ldots,h_n^{-1}-1$, replaced by the estimators~\eqref{pilot} and~\eqref{noiselevelest}, respectively.
\end{corro}
Unlike in \eqref{cltspoteq}, in which we obtain a mixed normal limiting distribution, the matrix $\mathcal{Z}$ is completely known. It is given by twice the ``symmetrizer matrix'' introduced by \citet[ch. 11]{magnus} and corresponds to the covariance structure of the empirical covariance of a $d$-dimensional (standard) Gaussian vector.

The asymptotic variance-covariance matrix in~\eqref{cltspoteq} is the same instantaneous process that appears integrated over $[0,1]$ as variance-covariance matrix of the integrated covariance matrix estimator in \cite{BHMR}. Accordingly, Theorem~\ref{cltspot} is in line with the results on classical realized volatility in the absence of noise for $d=1$ and the nonparametric Nadaraya-Watson-type kernel estimator by \cite{kristensen} with asymptotic variance $2\sigma_s^4\int_{\mathds{R}}k^2(z)\,dz$, where $k$ denotes the used kernel. In our case, the estimator is of histogram-type and the rectangle kernel does not appear in the asymptotic variance. Let us point out that estimator \eqref{spotcov}, building on optimal combinations over spectral frequencies, is more advanced than a usual histogram-estimator. When comparing our nonparametric estimator~\eqref{spotcov}, e.g., to the aforementioned estimator by \cite{kristensen}, in our case, the actual bandwidth is $(2K_n+1) h_n$ (or smaller), since we smooth over (up to) $(2K_n+1)$ adjacent blocks of length $h_n$. In this context,  one can as well think of employing $h_n^{-1}$ de-noised block statistics as underlying observations.

Regarding the convergence rate in~\eqref{cltspoteq}, we may focus on the case $\alpha= 1/2$, which is tantamount to the spot volatility matrix process $(\sigma_s)_{s\in[0,1]}$ being as smooth as a continuous semi-martingale. This assumption yields the rate $n^{1/8-\varepsilon}$, for any $\varepsilon>0$, such that we almost attain the optimal rate $n^{1/8}$, which is obviously lower than the corresponding rate for \textit{integrated } (co-)variance estimators in the setting with noise, $n^{1/4}$ \citep[see][]{hoff2012}. Notably, our spot covariance matrix estimator~\eqref{spotcov} converges considerably faster than existing noise-robust spot volatility estimators based on the difference quotient of integrated volatility  estimates \citep[e.g.][]{zu14}. The two-step approach \eqref{spotcov} with combinations over different frequencies strongly reduces the estimator's variance (compared to simpler methods). This is well confirmed in our finite-sample simulations in Section \ref{sec:sim}.

Theorem~\ref{cltspot} and Corollary~\ref{fclt} hold for estimation points $s\in[0,1]$, both in the interior and in the boundary region of the unit interval.  This result is a consequence of the estimators~\eqref{pilot} and~\eqref{spotcov} being of histogram-type, implying that smoothing is conducted by averaging over a set of adjacent blocks. The latter merely needs to contain time $t$, and does not have to be centered around the point of estimation. 

Finally, Theorem~\ref{cltspot} may be employed to deduce asymptotic results for the estimators of spot correlations and spot betas. These can be considered as the instantaneous counterparts to the integrated quantities studied, e.g., in \citet{abdl2003} and \citet{barnd2004}. In this context, focus on those elements of the spot covariance matrix $\Sigma_t,t\in[0,1]$, involving only the indices $p,q\in\left\{1,\ldots,d\right\}$. Further, denote the spot correlation and beta estimators based on~\eqref{spotcov} by $\hat\rho_s^{(pq)}=\hat\Sigma_s^{(pq)}/\sqrt{\hat\Sigma_s^{(pp)}\hat\Sigma_s^{(qq)}}$ and $\hat\beta_s^{(pq)}=\hat\Sigma_s^{(pq)} / \hat\Sigma_s^{(pp)}$. Then, Theorem~\ref{cltspot} implies by application of the Delta-method that
\begin{subequations}
\begin{align}
\label{cltspotcorr}
n^{\beta/2}\big(\hat\rho_s^{(pq)}-\rho_s^{(pq)}\big)\stackrel{d-(st)}{\longrightarrow}\mathbf{N}\Big(0,\mathds{AV}_{\rho,s}^{(pq)}\Big),\quad s\in[0,1]\,,\\
\label{cltspotbeta}
n^{\beta/2}\big(\hat\beta_s^{(pq)}-\beta_s^{(pq)}\big)\stackrel{d-(st)}{\longrightarrow}\mathbf{N}\Big(0,\mathds{AV}_{\beta,s}^{(pq)}\Big),\quad s\in[0,1]\,,
\end{align}
with
\begin{align}
\mathds{AV}_{\rho,s}^{(pq)}&=\Sigma_s^{(pp)}\Sigma_s^{(qq)}\mathds{AV}_s^{(p-1)d+q,(p-1)d+q}+\frac{\big(\Sigma_s^{(pq)}\big)^2}{4\big(\Sigma_s^{(pp)}\big)^3\Sigma_s^{(qq)}} \,\mathds{AV}_s^{(p-1)d+p,(p-1)d+p}\\\nonumber
&\quad +\frac{\big(\Sigma_s^{(pq)}\big)^2}{4\big(\Sigma_s^{(qq)}\big)^3\Sigma_s^{(pp)}} \,\mathds{AV}_s^{(q-1)d+q,(q-1)d+q}-\frac{\Sigma_s^{(pq)}}{\Sigma_s^{(pp)}} \,\mathds{AV}_s^{(p-1)d+q,(p-1)d+p}\\\nonumber
&\quad -\frac{\Sigma_s^{(pq)}}{\Sigma_s^{(qq)}} \,\mathds{AV}_s^{(p-1)d+q,(q-1)d+q}+\frac{\big(\Sigma_s^{(pq)}\big)^2}{2\big(\Sigma_s^{(pp)}\Sigma_s^{(qq)}\big)^2} \,\mathds{AV}_s^{(p-1)d+p,(q-1)d+q}\,,
\\
\mathds{AV}_{\beta,s}^{(pq)}&=\big(\Sigma_s^{(pp)}\big)^{-2}\mathds{AV}_s^{(p-1)d+q,(p-1)d+q}+\big(\Sigma_s^{(pq)}\big)^2\big(\Sigma_s^{(pp)}\big)^{-4}\mathds{AV}_s^{(p-1)d+p,(p-1)d+p}\\\nonumber
 &\quad -2\Sigma_s^{(p,q)}\big(\Sigma_s^{(p,p)}\big)^{-3}\mathds{AV}_s^{(p-1)d+q,(p-1)d+p}\,,
\end{align}
\end{subequations}
where $\mathds{AV}_s$ denotes the asymptotic variance-covariance matrix in~\eqref{cltspoteq}. Feasible versions of the central limit theorems~\eqref{cltspotcorr} and~\eqref{cltspotbeta} can be readily obtained analogously to Corollary~\ref{fclt}.

\subsection{Choice of Inputs\label{sec:lmm_input}}

The proposed spot covariance matrix estimator~\eqref{spotcov} depends on  four  input parameters to be chosen: (i)~the block length $h_n$, (ii)~the maximum spectral frequency $J_n$, (iii)~the maximum frequency for the pre-estimator~\eqref{pilot}, $J_{n}^p$, as well as (iv)~the length of the smoothing window,  $K_n$.

For (i) , Theorem~\ref{cltspot} requires that $h_n=\mathcal{O}\big(\log{(n)}n^{-1/2}\big)$. (ii) is given by $\lfloor \min_p n_p h_n\rfloor $, but a spectral cut-off $J_n=\mathcal{O}\big(\log{(n)}\big)$ can be chosen, since the optimal weights decay fast with increasing frequency $j$, making higher frequencies asymptotically negligible. The effect of quickly diminishing optimal weights implies that (iii) should be fixed at a value not ``too large'', e.g., $J_{n}^{p}=5$. The reason is that the cut-off directly determines the (uniform) weights in the pre-estimator~\eqref{pilot}. For (iv), we generally set $K_n=\KLEINO\big( n^{\alpha/(2\alpha+1)}\big)$. The latter choice implies undersmoothing, thereby forfeiting rate-optimality of the estimator, but provides us a central limit theorem.
Under the ``continuous semi-martingale or smoother'' assumption ($\alpha\geq 1/2$) for the spot volatility matrix process, which seems admissible in most financial applications, we set $K_n=\mathcal{O}\big(n^{1/4-\varepsilon}\big)$ for some $\varepsilon>0$.

In practice, we introduce proportionality parameters for (i), (ii) and (iv), i.e. $h_n=\theta_h\log{(n)}n^{-1/2}$, $J_n=\lfloor\theta_J \log{(n)}\rfloor$ and $K_n=\lceil\theta_K n^{1/4-\delta}\rceil$, where $\theta_h,\theta_J,\theta_K>0$ and $\delta$ denotes a small positive number. 
We discuss the specific choice of the above input parameters in more detail in Sections~\ref{sec:sim} and \ref{sec:empspot}.

\section{Simulation Study\label{sec:sim}}

We conduct a simulation study to examine the following issues. First, we analyze the impact of different choices of the input parameters $\theta_h$, $\theta_J$ and $\theta_K$ introduced in Section~\ref{sec:lmm_input} on the estimator's finite-sample performance. We consider different scenarios which mimic both ``regular'' trading days as well as ``unusual'' trading days in periods of financial stress. Second, we investigate the frequency of non-positive semi-definite estimates and whether simple eigenvalue truncation techniques translate into an improved finite-sample precision. 

We consider a high-dimensional setting with $d=15$. For 15 assets, we estimate a 120-dimensional volatility matrix and the estimator utilizes weight matrices with 7260 entries. To ensure parsimony in this framework, we assume that the efficient log-price process follows a simple factor structure as employed, e.g., in \citet{bn1}. We extend the latter to incorporate both a flexible stochastic and a non-stochastic seasonal volatility component, which is modeled by a Flexible Fourier Form as introduced by \citet{gall81}. We dilute the observations of the efficient log-price process by serially dependent microstructure noise with $R=1$. Finally, asynchronicity effects are introduced by drawing the observation times $t_i^{(p)}$, $i=1,\ldots,n_p$, from independent Poisson processes. Details on the simulation setting are provided in the web appendix, Section 2.

To investigate the impact  of the chosen input parameters, we compute the LMM estimator~(\ref{spotcov}) over a grid of values for $\theta_h$, $\theta_J$ and $\theta_K$. For each combination and in each replication $m=1,\ldots,M$, we evaluate the (normalized) mean integrated Frobenius distance between the resulting estimates $\hat\Sigma_{t,m}$  and their ``true'' counterparts.  Hence, we compute
\begin{align}\label{mifb}
\text{MIFB}\defeq \left(M d^2\right)^{-1}\sum_{m=1}^M \int_0^1\sum_{p,q=1}^d\left[\hat\Sigma^{(pq)}_{t,m}/\Sigma^{(pq)}_{t,m}-1\right]^2 dt, 
\end{align}
where $M$ is the number of replications.  In addition, we evaluate the average normalized mean integrated squared errors of the variance and covariance estimates, respectively, i.e., 
\begin{align}\label{misec}
\text{MISE}_{c}&\defeq \left(M d (d-1)/2\right)^{-1}\sum_{m=1}^M \int_0^1\sum_{p\neq q}\left[\hat\Sigma^{(pq)}_{t,m}/\Sigma^{(pq)}_{t,m}-1\right]^2 dt,\\\label{misev} 
\text{MISE}_{v}&\defeq \left(M d\right)^{-1}\sum_{m=1}^M \int_0^1\sum_{p=1}^d\left[\hat\Sigma^{(pp)}_{t,m}/\Sigma^{(pp)}_{t,m}-1\right]^2 dt.
\end{align}
Finally, to examine how often the estimator~(\ref{spotcov}) yields non-positive semi-definite estimates, we compute the percentage of replications in which all spot covariance matrix estimates are positive semi-definite.

Panels~A,~B and~C of Table~\ref{tab:sim_theta_choice} report the values  of the input parameters minimizing MIFB, MISE$_c$ and MISE$_v$, respectively, for $M=3000$ along with the square roots of the latter distance measures. Panel~A additionally provides the performance implied by more ``extreme'' choices of the input parameters and the percentage of positive semi-definite estimates. Panels~B and ~C also report MISE$_c$ and MISE$_v$ based on the optimal parameter values  with respect to MIFB. The MIFB-optimal values of the input parameters yield a configuration with (on average) $\lceil h_n^{-1} \rceil=80$ blocks spanning about 5 minutes each, a spectral cut-off  $J_n=19$ and a smoothing window of $K_n=6$ blocks. Regarding deviations from the MIFB-optimal values of the input parameters, considerable precision losses occur in only two cases. First, when setting $\theta_h$ extremely low, resulting in more than $480$ blocks per day on average. Second, for a very small choice of $\theta_J$, as spectral frequencies are cut off too early. In particular, we observe that the two-step method~\eqref{spotcov} clearly outperforms a simple histogram-type estimator, which relies only on the first frequency $J_n=1$.  We can conclude that the performance of the (full) spot covariance matrix estimator is quite robust for a range of sensible input choices.

\begin{table}
\caption{Performance of LMM spot covariance matrix estimator depending on $\theta_h$, $\theta_J$ and $\theta_K$. 
RMIFB, RMISE$_c$ and RMISE$_v$ are the square roots of~\eqref{mifb},~\eqref{misec} and~\eqref{misev}, respectively, computed based on $M=3000$ Monte Carlo replications and reported in percentage points. $\%$ PSD denotes the percentage of Monte Carlo replications yielding exclusively positive semi-definite spot covariance matrix estimates. Opt* corresponds to inputs, which are optimal with respect to RMIFB.}\label{tab:sim_theta_choice}
\centering
\small
\begin{tabular*}{\textwidth}{@{\extracolsep{\fill}}lrrrrr}
\toprule
& \multicolumn{5}{c}{Full Covariance Matrix}\\
\midrule
& $\theta_h$ & $\theta_J$ & $\theta_K$ & RMIFB&$\%$ PSD \\
\midrule
Opt & 0.150 & 6.000 & 2.000 & 24.410 & 78.300 \\
  & 0.150 & 1.000 & 2.000 & 35.627 & 99.900 \\
  & 0.150 & 10.000 & 2.000 & 24.444 & 77.300 \\
  & 0.150 & 6.000 & 1.200 & 24.910 & 86.400 \\
  & 0.150 & 6.000 & 4.800 & 26.835 & 57.600 \\
  & 0.025 & 6.000 & 2.000 & 47.043 & 0.000 \\
  & 0.250 & 6.000 & 2.000 & 25.331 & 98.500 \\
 \midrule
 & \multicolumn{5}{c}{Covariances}\\
\midrule
& $\theta_h$ & $\theta_J$ & $\theta_K$ & $\text{RMISE}_{c}$& \\
\midrule
Opt & 0.100 & 7.000 & 2.400 & 22.944 &\\
Opt* & 0.150 & 6.000 & 2.000 & 23.050 &\\
 \midrule
 & \multicolumn{5}{c}{Variances}\\
\midrule
& $\theta_h$ & $\theta_J$ & $\theta_K$ & $\text{RMISE}_{v}$& \\
\midrule
Opt & 0.150 & 2.000 & 1.200 & 17.720 \\
Opt* & 0.150 & 6.000 & 2.000 & 22.846 \\
\bottomrule
\end{tabular*}
\end{table}

When focusing on covariance estimates only, MISE$_c$-optimal values of the input parameters would mainly imply an increase in the average number of blocks to around $120$ per day and a corresponding lengthening of the smoothing window to $7$ blocks,  while precision remains very close to the one implied by MIFB-optimal inputs. For the variances, the spectral cutoff would reduce to around $6$, while the smoothing window would shorten to around $4$ blocks. The MIFB-optimal values imply a non-negligible increase in  MISE$_v$. Table~\ref{tab:sim_theta_choice} further shows that employing MIFB-optimal values of the input parameters yields positive semi-definite spot covariance matrix estimates in around  $78\%$ of the cases, while increasing the spectral cutoff or reducing the block length leads to more cases with non-positive semi-definite outcomes.

Simulation results based on a LMM estimator with a truncation of negative eigenvalues at zero, confirming an improved finite-sample performance, and for a setting with volatility (co-)jumps can be found in the web appendix.

\section{Empirical Study\label{sec:emp}}
\subsection{Implementation\label{sec:impl}}
We apply the estimators presented in Section \ref{sec:estimspotcov} to our dataset described in Section \ref{sec:data}. In order to obtain spot covariance matrix estimates for the entire trading day including the period immediately after the start of trading, we initially consider the two-sided version of the estimator  \eqref{spotcov}. We then use the mid-quote revisions for the Nasdaq 100 constituents to estimate $30 \times 30$ spot covariance matrices, yielding pair-wise spot covariances and correlations, as well as individual volatilities. We select the relevant inputs as discussed in Section~\ref{sec:lmm_input}. The corresponding proportionality parameters are set to the values found to be ``optimal'' in the extended simulation study given in Section 2.2 of the web appendix. Hence, we set  $\theta_h=0.175$, $\theta_J=7$, $\theta_K=2$ and $J_n^p=5$.

Section 4 of the web appendix reports summary statistics for the number of blocks, the spectral cut-off and the length of the smoothing window as induced by the underlying data for both the entire sample and each year. On average, we use approximately $27$ blocks per day, resulting in an average block length of $14$ minutes. Spectral frequencies are cut off at nearly $48$, while the average length of the smoothing window is about $6$ blocks, translating into roughly 80 minutes.
\subsection{Intraday Behavior of Spot (Co-)Variances\label{sec:empspot}}
As a first step, we investigate the presence of seasonality effects in spot (co-)volatilities. Seasonal patterns in intraday volatilities have been confirmed, e.g., in the seminal studies by \citet{ande97,ande98}. For equity returns, volatilities typically exhibit a U-shape, i.e., volatility is higher at the opening and before the closure of the market, while being lower around midday. Similar effects have been documented for other measures of intraday trading activity such as bid-ask spreads \citep[e.g.][]{chan1995}, durations between trade and quote arrivals \citep[e.g.][]{engle1998} and transaction volumes \citep[e.g.][]{brown2011}.

\begin{figure}
\centering
\hspace*{-4ex}\hfill\subfigure[Covariances]{\includegraphics[width=0.42\textwidth]{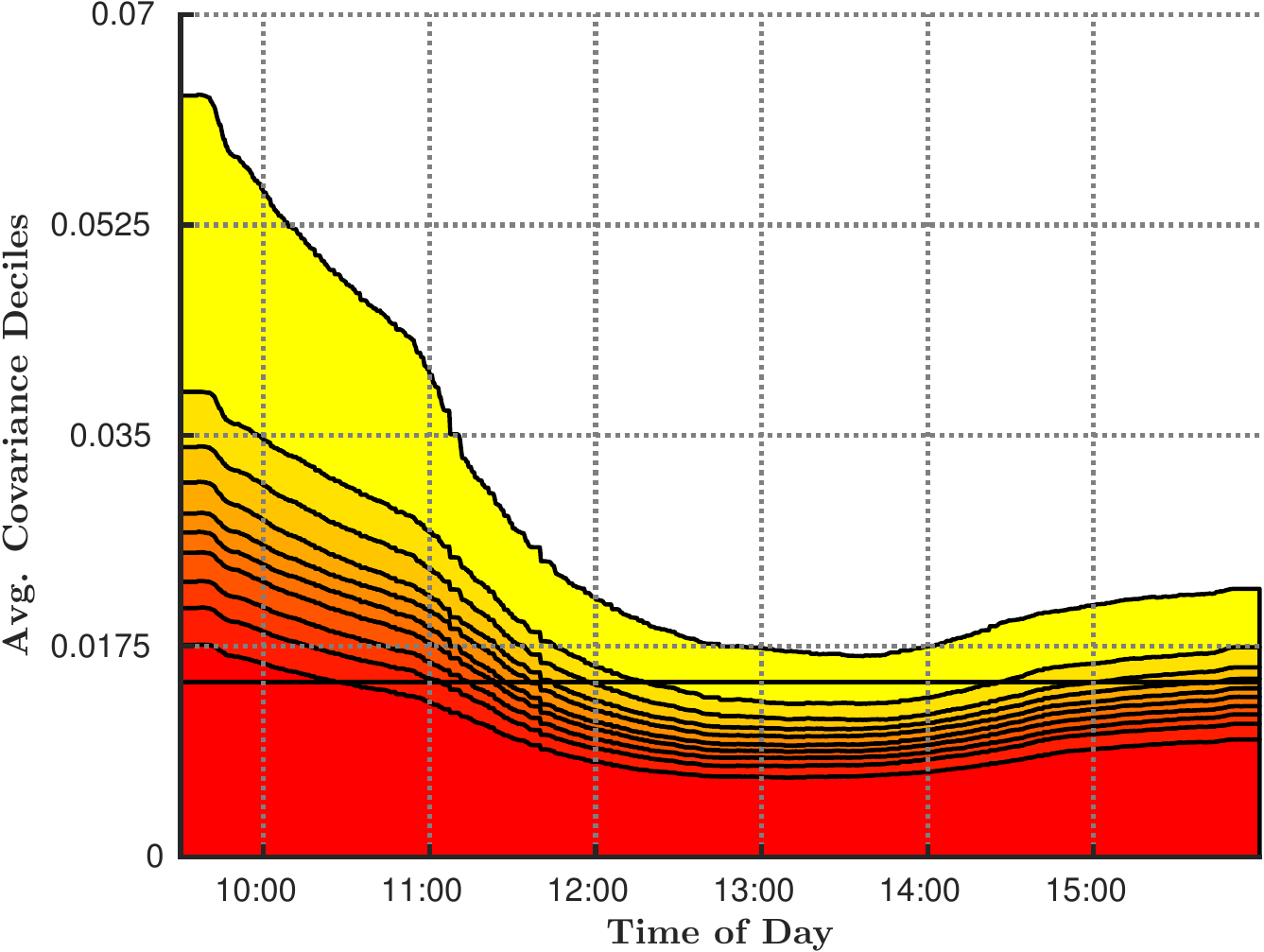}}\hspace*{-0.3ex}\hfill
\subfigure[Correlations]{\includegraphics[width=0.42\textwidth]{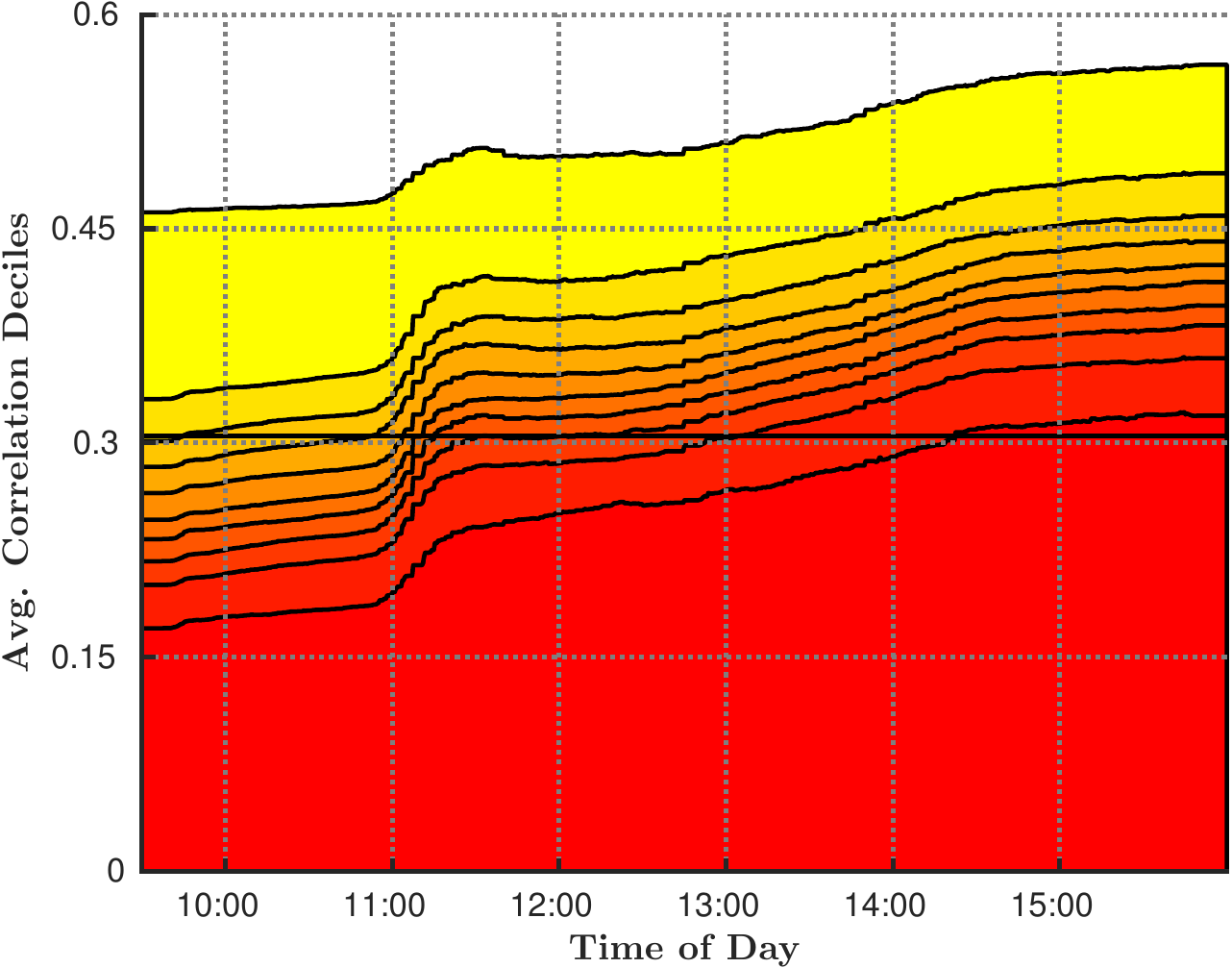}}\hfill\\
\subfigure[Volatilities]{\includegraphics[width=0.42\textwidth]{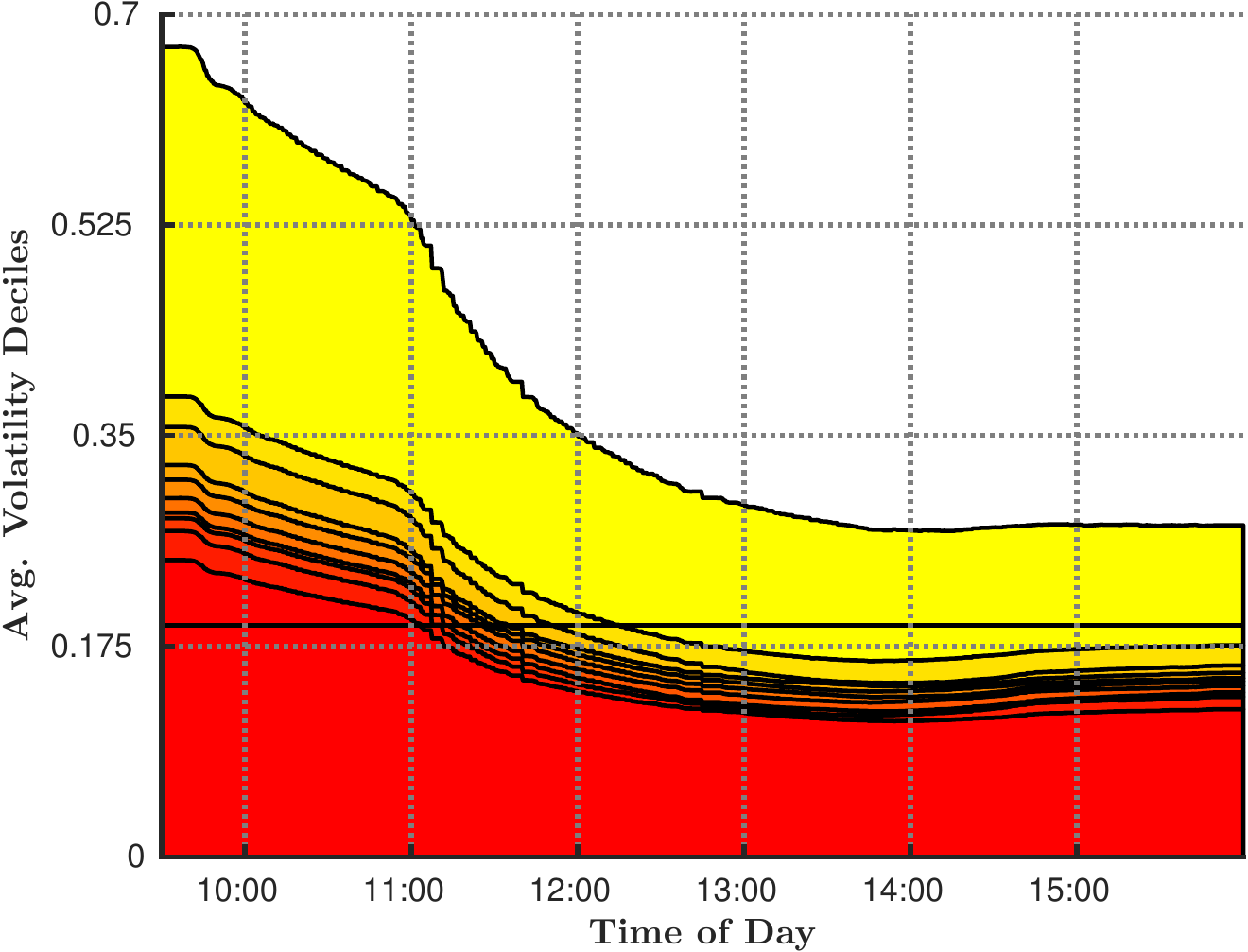}}\hfill
\caption{Cross-sectional deciles of across-day averages of spot covariances, correlations and volatilities. Spot estimates are first averaged across days for each asset pair. Subsequently, cross-sectional sample deciles of the across-day averages are computed. The solid horizontal line corresponds to the cross-sectional median of the across-day averages of \textit{integrated} covariance, correlation and volatility estimates. These are based on the LMM estimator of the integrated (open-to-close) covariance matrix by \citet{BHMR} accounting for serially dependent noise and using the same input parameter configuration as the spot estimators. ``Unusual days'' discussed in Section~\ref{sec:event} as well as days with scheduled FOMC announcements are removed. Covariances and volatilities are annualized.}
\label{fig:cov_corr_vol_fulls}
\end{figure}

Figure \ref{fig:cov_corr_vol_fulls} shows the cross-sectional deciles of across-day averages of spot covariances and correlations for each asset pair as well as volatilities for each asset. The averages were computed while omitting the ``unusual'' days analyzed in Section~\ref{sec:event} as well as days with scheduled announcements of the federal funds rate target by the Federal Open Market Committee (FOMC).  A more detailed discussion of the latter can be found below.

We observe distinct intraday seasonality patterns, which, interestingly, do not only apply to volatilities. Rather, covariances clearly decline at the beginning of the trading day, stabilize around noon on a widely constant level and slightly increase before market closure. Interestingly, the resulting correlations show a reverse pattern and significantly increase during the first trading hour. The latter is caused by spot volatilities that decay faster than the corresponding covariances at the beginning of the trading day. Hence, the (co-)variability between assets is highest after start of trading which might be caused by the processing of common information analogously to the higher overall inflow of public and private information during that period \citep[see, e.g.,][]{hasb1991,madh1997}. The latter effect, however, appears to imply an even more pronounced increase in assets' idiosyncratic risk as reflected by spot volatilities,  overcompensating  higher covariances and leading to lower correlations at the beginning of the trading day. Interestingly, spot volatilities drop significantly faster than underlying covariances during the first trading hour. Shortly after opening, spot volatilities are approximately twice as high as the (average) \emph{daily} volatility (computed based on the open-to-close integrated variance estimate), but strongly decline thereafter. This makes correlations sharply increasing between 10:00 and 11:00 am. Accordingly, we observe that median spot correlations range between approximately $0.2$ and $0.4$ across a day. This is in contrast to a \emph{daily} correlation (computed from the open-to-close integrated covariance estimate) of approximately $0.3$ and shows that even on average, intraday variability of correlations and covariances is substantial. Finally, we repeat the above analysis on a year-by-year basis. The results reported in Section 5 of the web appendix show that the intraday patterns mainly change in terms of level shifts over the years.

We can summarize the findings above as follows. First, spot covariances exhibit an intraday seasonality pattern closely resembling the U-shape which is typical for volatilities. Second, the combined diurnal patterns of spot covariances and volatilities imply that spot correlations tend to increase throughout the trading day. 

\begin{figure}
\centering
\hspace*{-4ex}\hfill\subfigure[Covariances]{\includegraphics[width=0.42\textwidth]{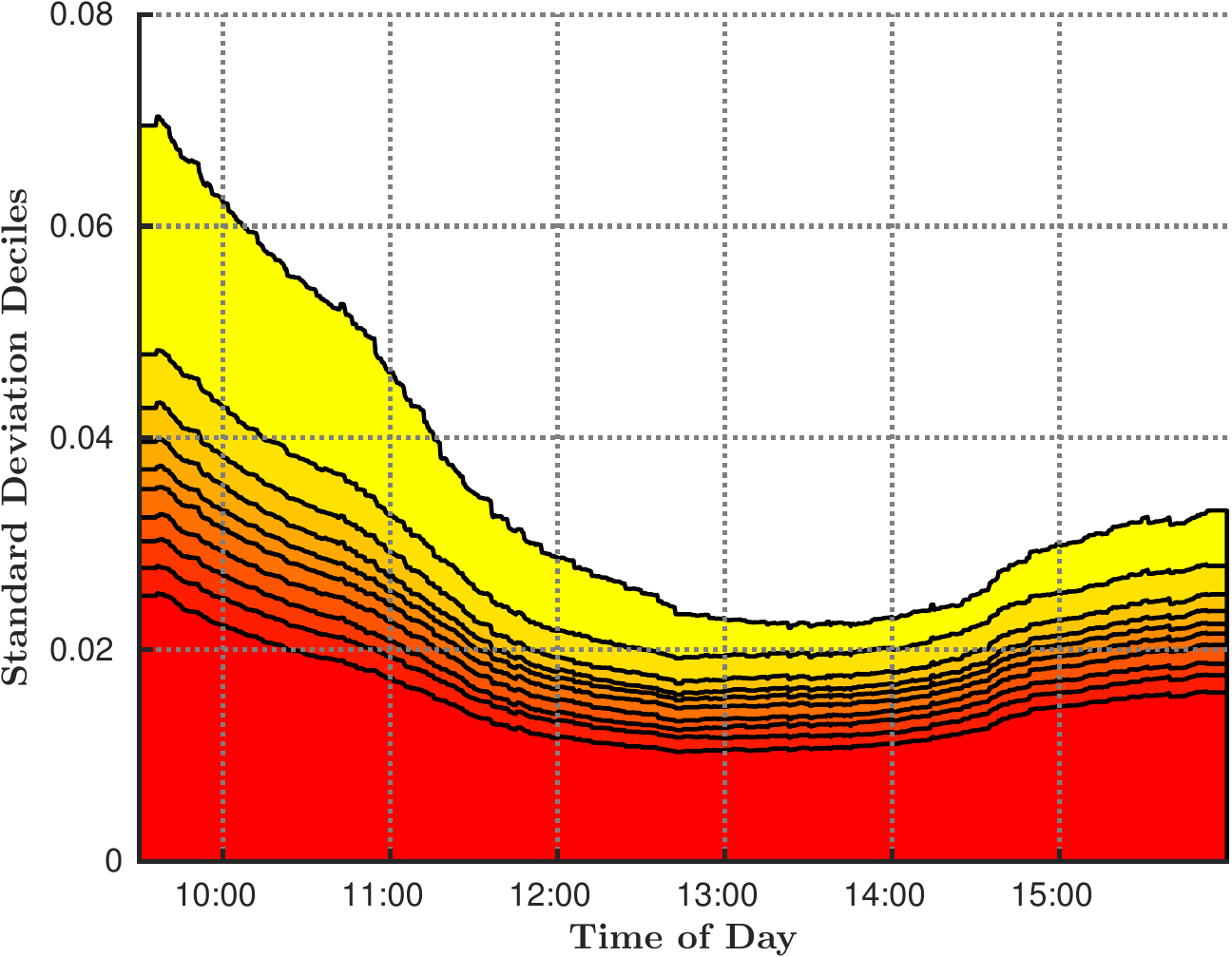}}\hspace*{-0.3ex}\hfill
\subfigure[Correlations]{\includegraphics[width=0.42\textwidth]{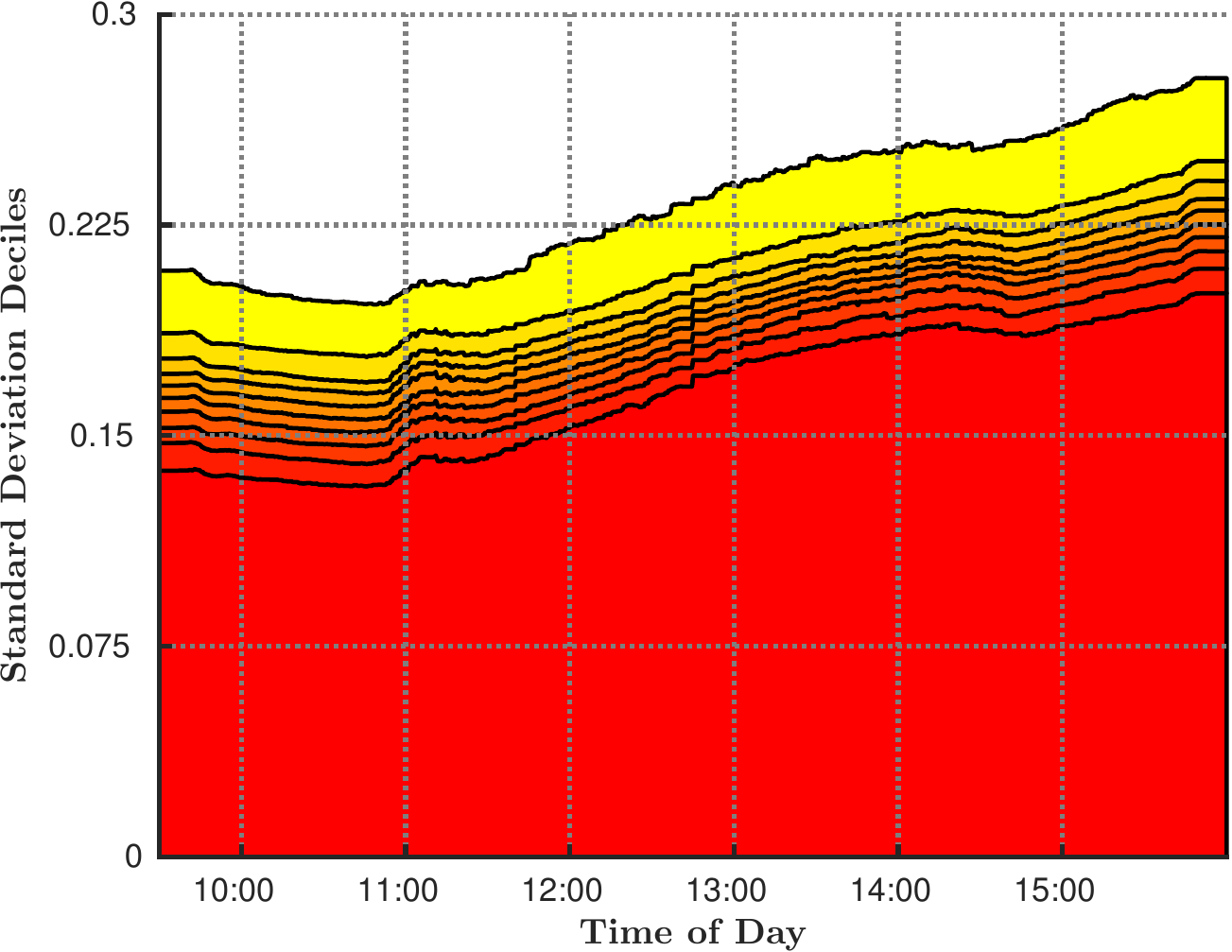}}\hfill\\
\subfigure[Volatilities]{\includegraphics[width=0.42\textwidth]{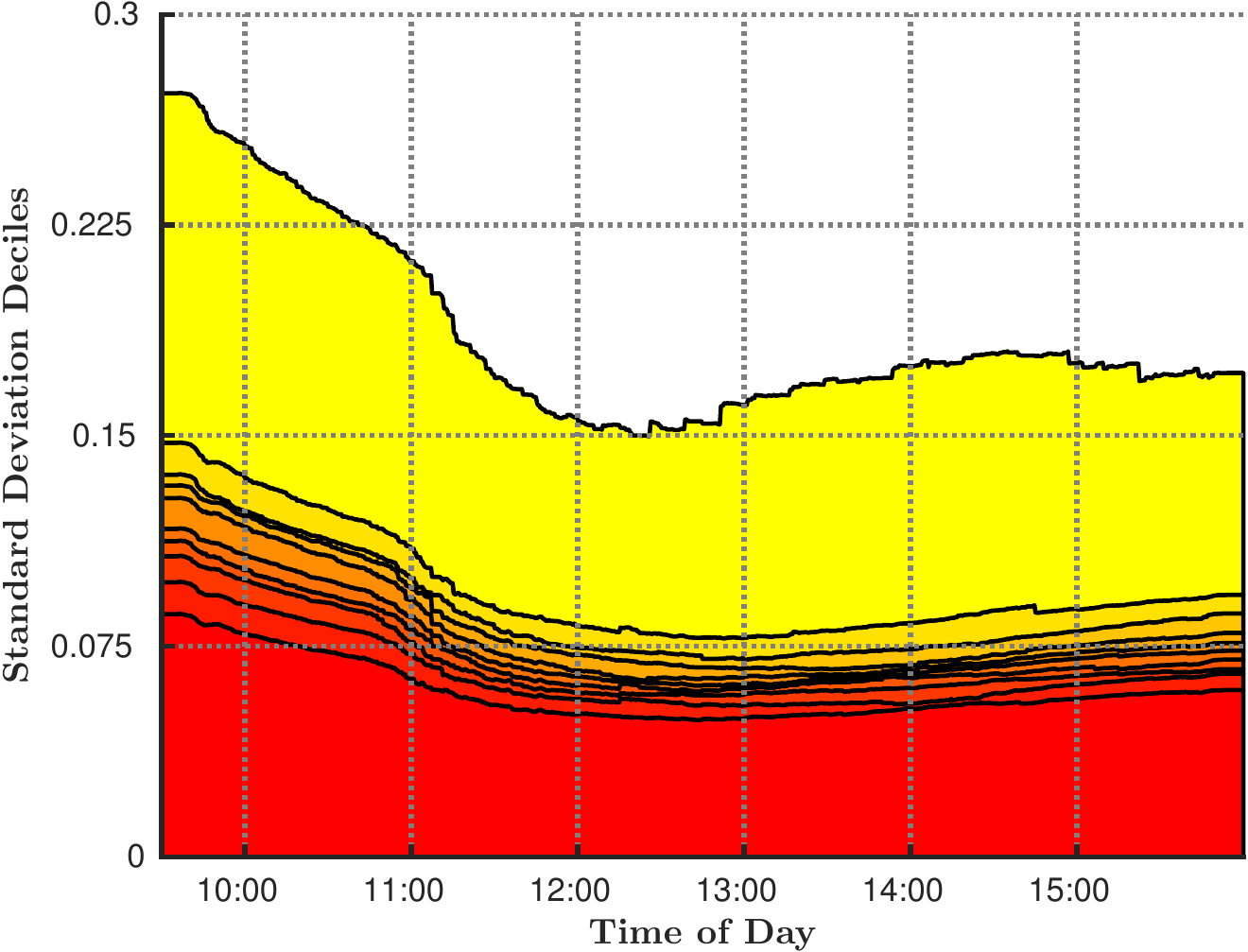}}\hfill
\caption{Cross-sectional deciles of across-day standard deviations of spot covariances, correlations and volatilities. First, sample standard deviations of spot estimates are computed across days for each asset pair. Subsequently, cross-sectional sample deciles of the across-day standard deviations are computed.  ``Unusual days'' discussed in Section~\ref{sec:event} as well as days with scheduled FOMC announcements are removed. Covariances and volatilities are annualized.}
\label{fig:cov_corr_vol_std}
\end{figure}

In Figure \ref{fig:cov_corr_vol_std}, we additionally compute, for each asset pair or asset and each point during the day, the standard deviation of spot covariances, correlations and volatilities \emph{across} days. We observe that the across-day variability in covariances is highest after market opening and shortly before closure. A similar picture is also observed for spot volatilities.  We associate the patterns  described above with effects arising from (overnight) information processed in the morning and increased trading activities in the afternoon, where traders tend to re-balance or close positions before the end of trading. Hence, idiosyncratic effects seem to become stronger during these periods, increasing the variability of (co-)variances. Interestingly, the across-day standard deviations in correlations show a reverse pattern. Thus, across-day variability in intraday correlations is lowest at the beginning of trading, increases until mid-day and is widely constant during the afternoon hours. Here, increased across-day covariance and volatility risk seem to compensate each other.

\begin{figure}
\centering
\hspace*{-4ex}\hfill\subfigure[Covariances]{\includegraphics[width=0.42\textwidth]{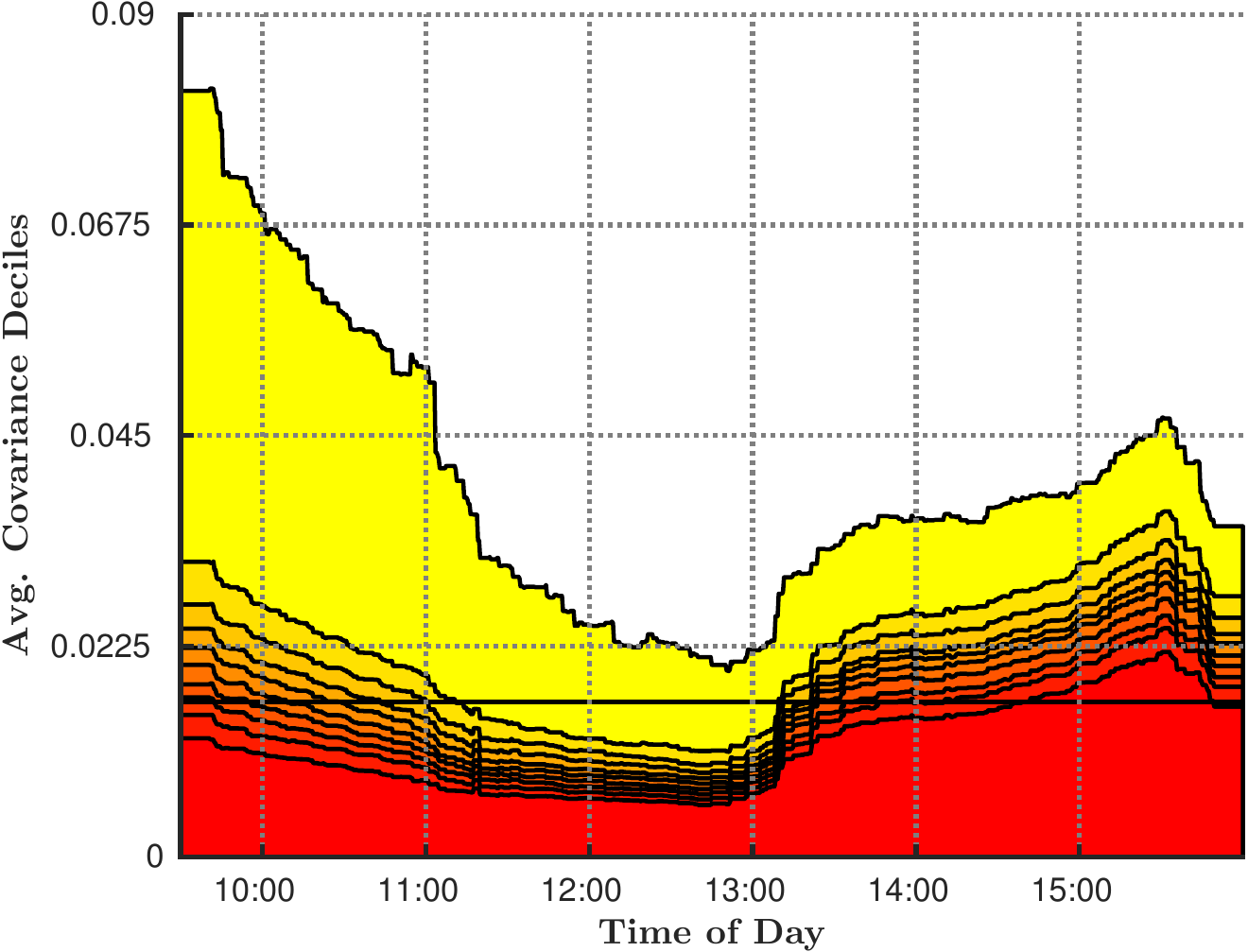}}\hspace*{-0.3ex}\hfill
\subfigure[Correlations]{\includegraphics[width=0.42\textwidth]{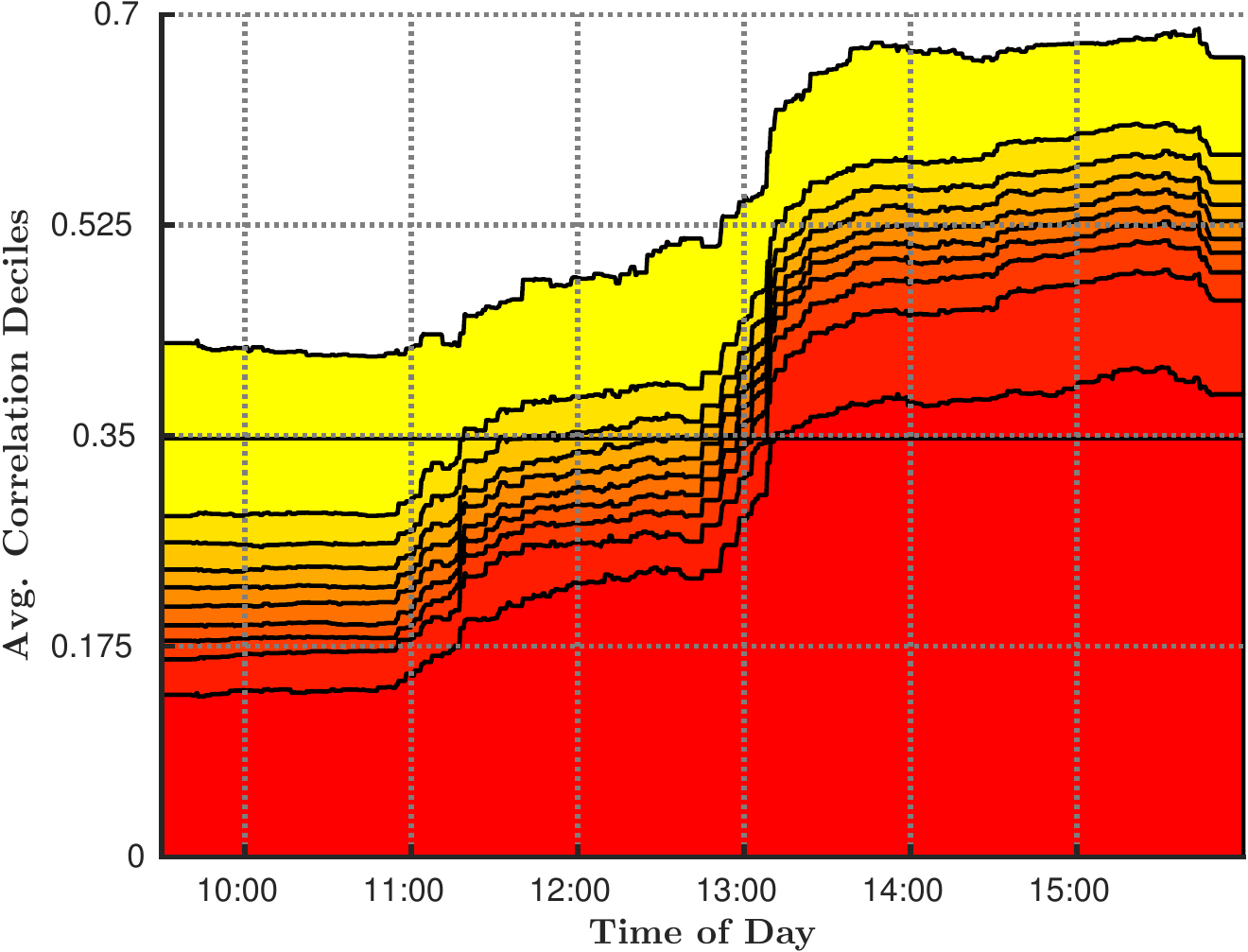}}\hfill\\
\subfigure[Volatilities]{\includegraphics[width=0.42\textwidth]{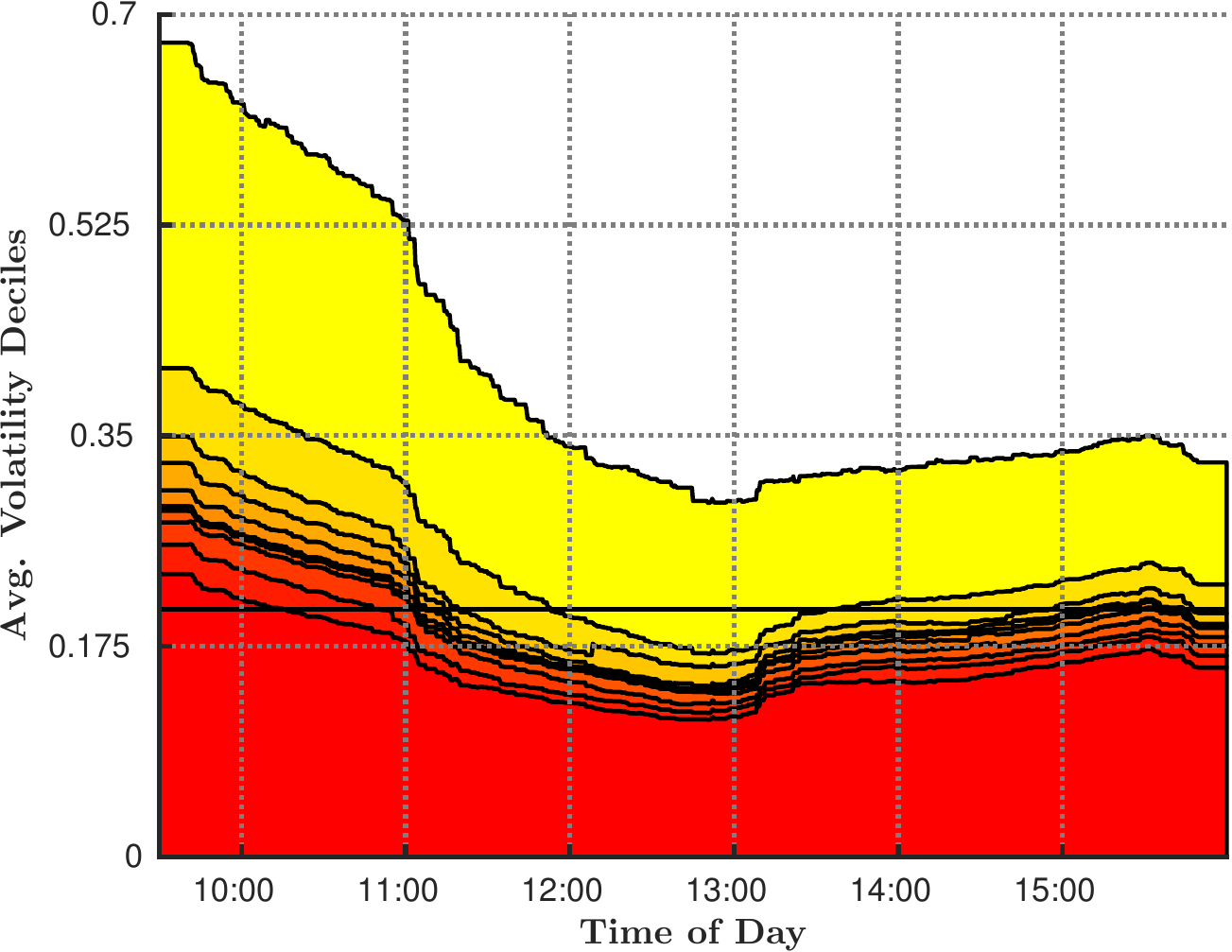}}\hfill
\caption{Cross-sectional deciles of across-day averages of spot covariances, correlations and volatilities on days with scheduled FOMC announcements. Spot estimates are first averaged across days for each asset pair. Subsequently, cross-sectional sample deciles of the across-day averages are computed. The solid horizontal line corresponds to the cross-sectional median of the across-day averages of \textit{integrated} covariance, correlation and volatility estimates. These are based on the LMM estimator of the integrated (open-to-close) covariance matrix by \citet{BHMR} accounting for serially dependent noise and using the same input parameter configuration as the spot estimators. Covariances and volatilities are annualized.}
\label{fig:cov_corr_vol_FOMC}
\end{figure}

It is well-known that equity returns respond to macroeconomic announcements both in terms of conditional means and volatilities \citep[see, e.g.,][]{ande2007,lund2009}. Accordingly, we compute the across-day averages analyzed above excluding days with major scheduled macroeconomic news announcements that regularly fall well within Nasdaq trading hours. For that purpose, we focuse on scheduled FOMC announcements, occurring at 2:15 pm  roughly every six weeks. 
For comparison, Figure~\ref{fig:cov_corr_vol_FOMC} reports the counterpart of Figure~\ref{fig:cov_corr_vol_fulls} based  on FOMC announcement days \emph{only}. Interestingly, we observe that around 1 pm, i.e. roughly one hour before the scheduled announcement, spot covariances exhibit a pronounced increase, while the rise in volatilities remains comparably modest. As a consequence, spot correlations simultaneously increase by a considerable extent.

\subsection{Two Unusual Days\label{sec:event}}
The previous section shows that spot correlations and covariances can substantially vary during a day, even if these patterns are averaged across time and assets. Here, we aim at analyzing the behavior of spot (co-)variability and the estimator~\eqref{spotcov}  in unusual market periods.  To prevent artifacts caused by ``forward-looking'' smoothing, we employ the one-sided version of the estimator. Further, we ensure a completely adaptive behavior of the latter by also estimating the long-run noise variance according to the method presented in the web appendix and determining the inputs following Section~\ref{sec:lmm_input} only based on observations available up to the point of estimation.

\begin{figure}
\centering
\hspace*{-4ex}\hfill\subfigure[Covariances]{\includegraphics[width=0.42\textwidth]{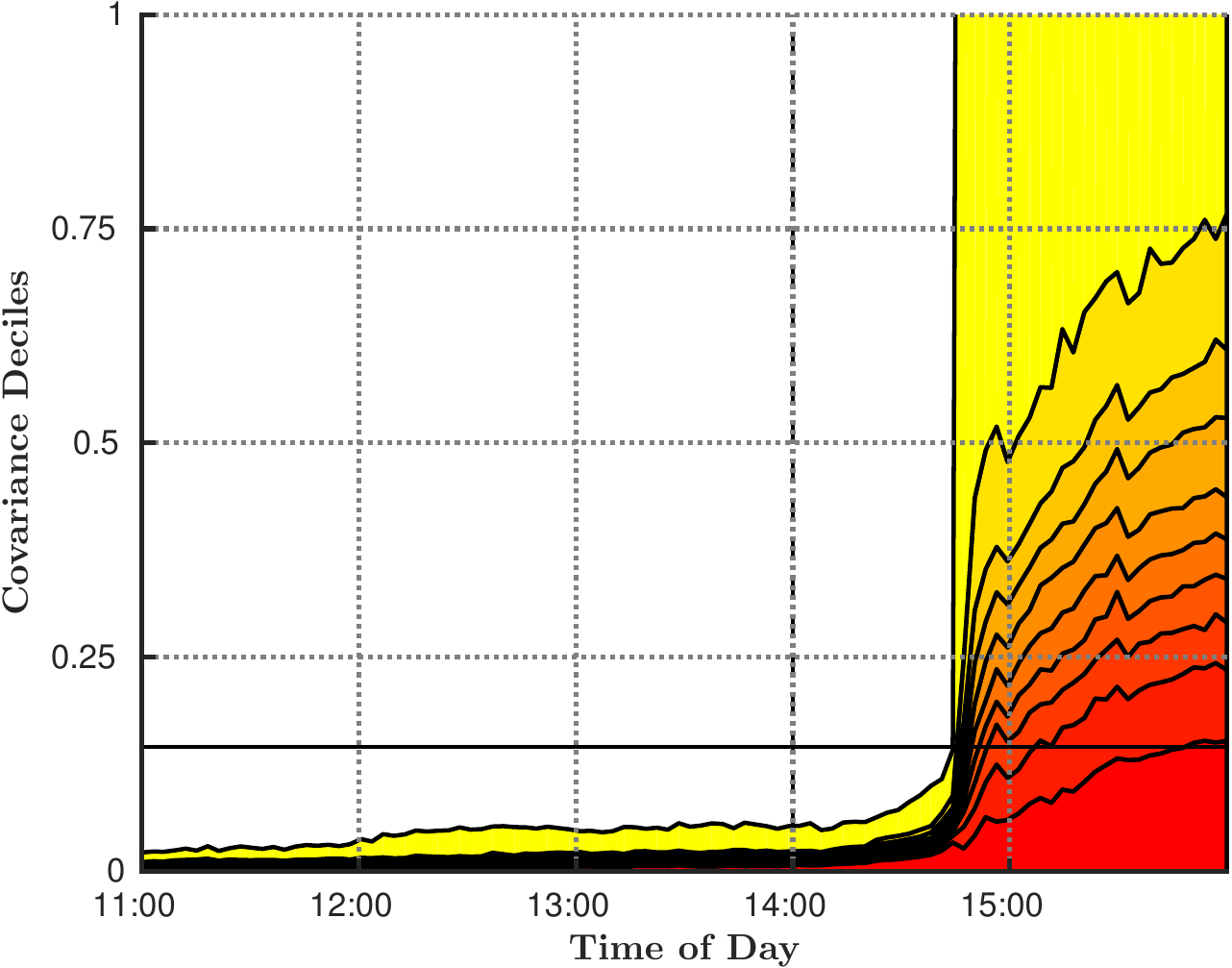}}\hspace*{-0.3ex}\hfill
\subfigure[Correlations]{\includegraphics[width=0.42\textwidth]{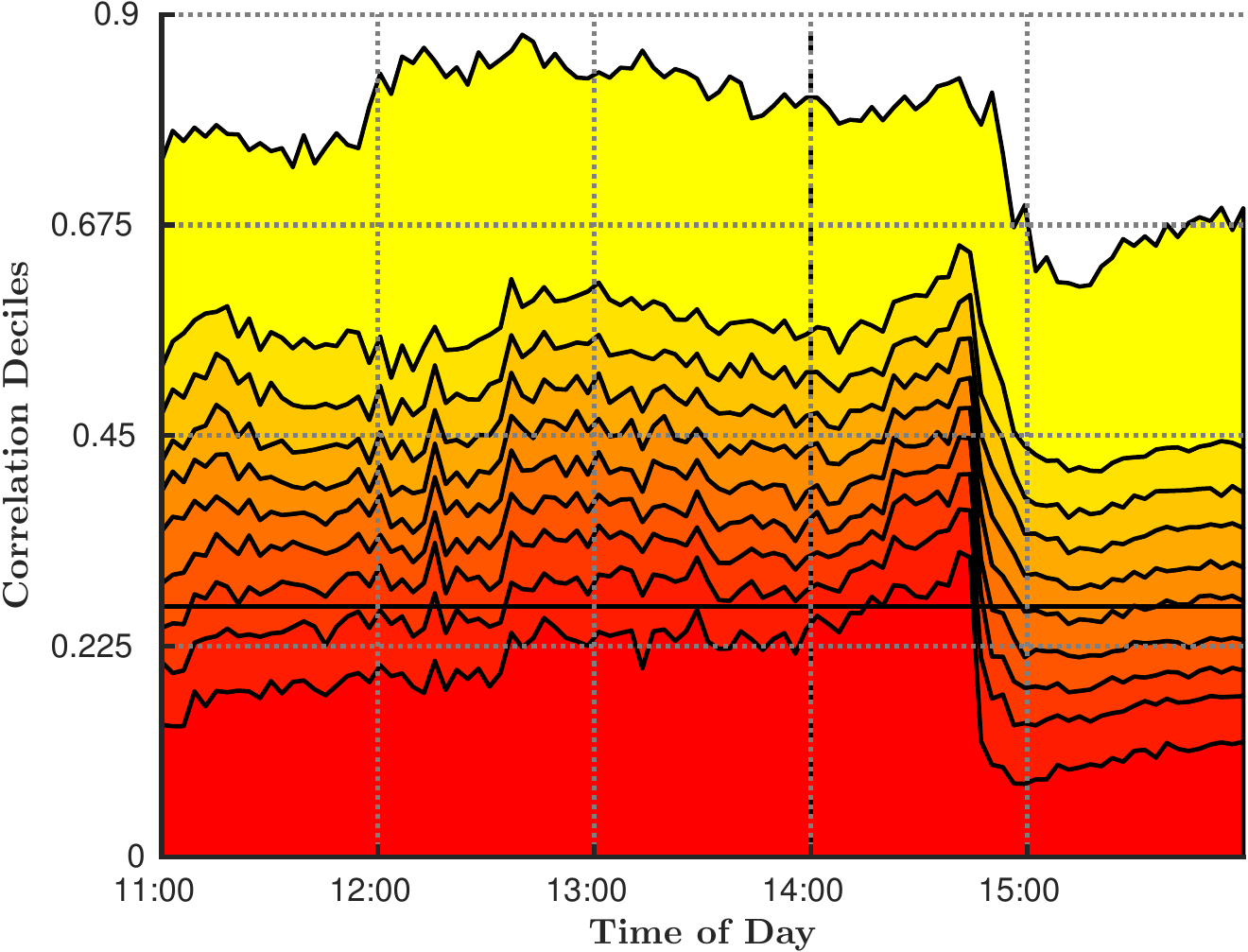}}\hfill\\
\subfigure[Volatilities]{\includegraphics[width=0.42\textwidth]{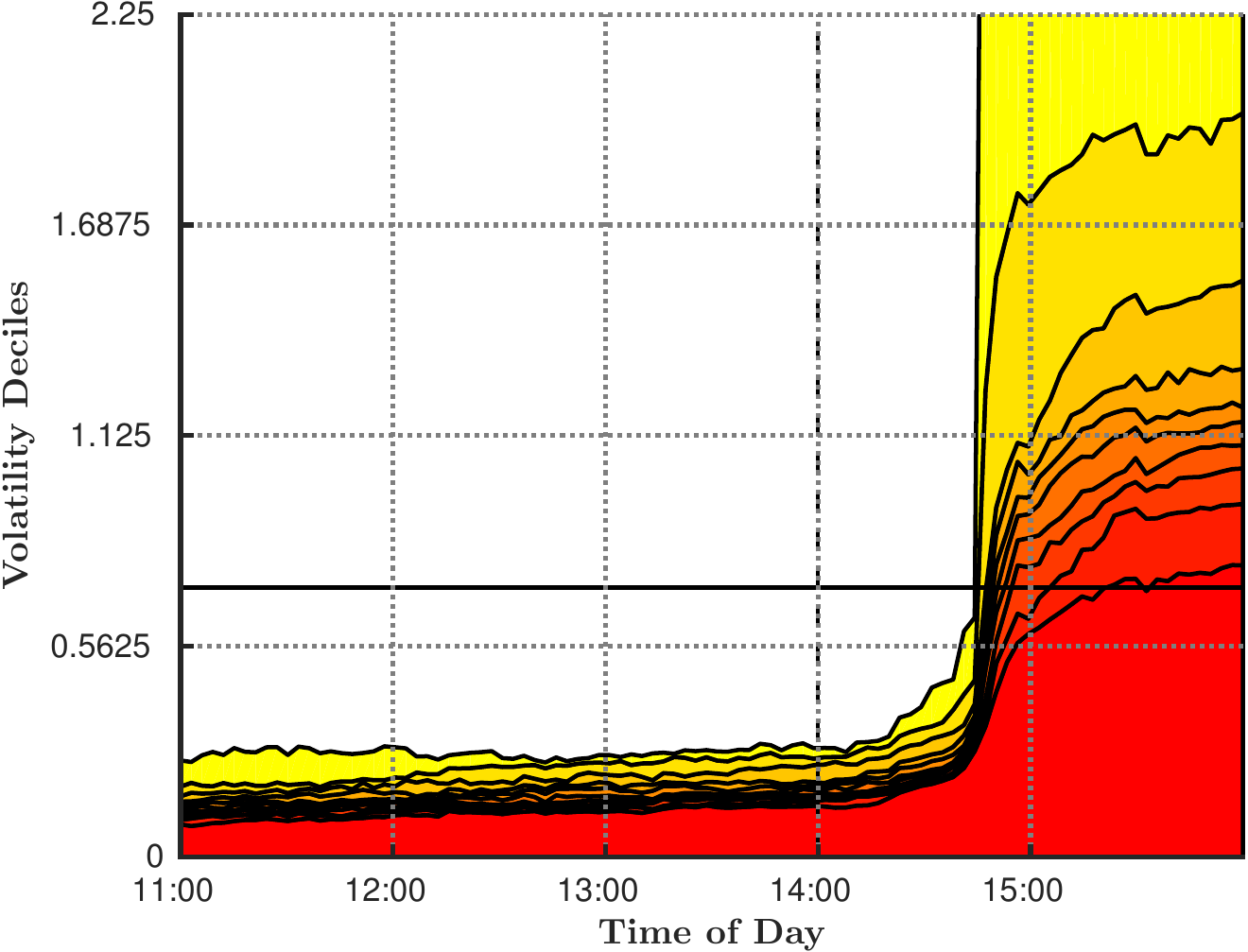}}\hfill
\caption{Cross-sectional deciles of spot covariances, correlations and volatilities (05/06/10). Solid horizontal line corresponds to the cross-sectional median of \textit{integrated} covariance, correlation and volatility estimates. These are based on the LMM estimator of the integrated (open-to-close) covariance matrix by \citet{BHMR} accounting for serially dependent noise and using the same input parameter configuration as the spot estimators.  Covariances and volatilities are annualized.}
\label{fig:cov_corr_vol_20100605}
\end{figure}

\begin{figure}
\centering
\hspace*{-4ex}\hfill\subfigure[Covariances/Correlations]{\includegraphics[height=0.34\textwidth]{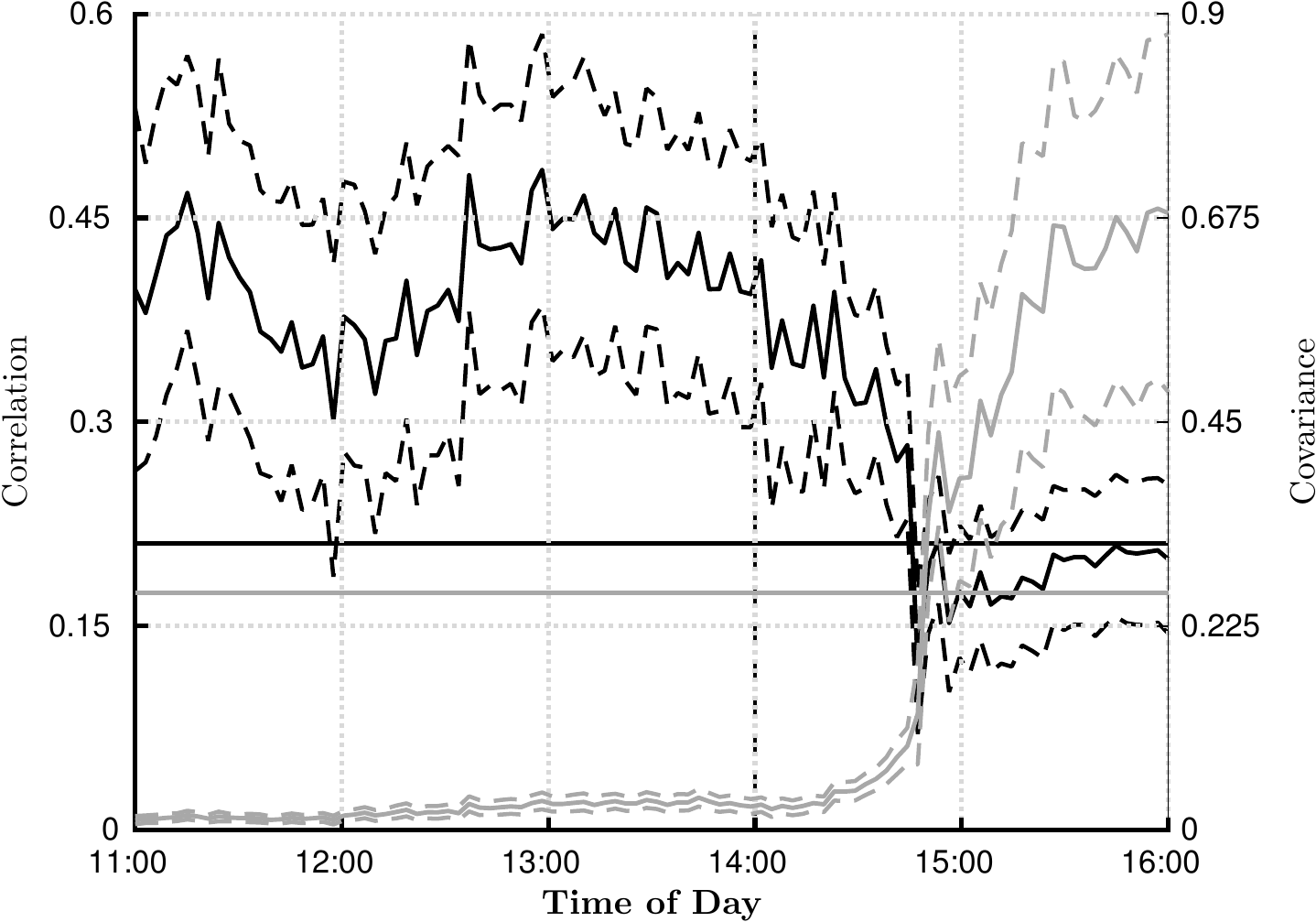}}\hfill
\subfigure[Volatilities]{\includegraphics[height=0.34\textwidth]{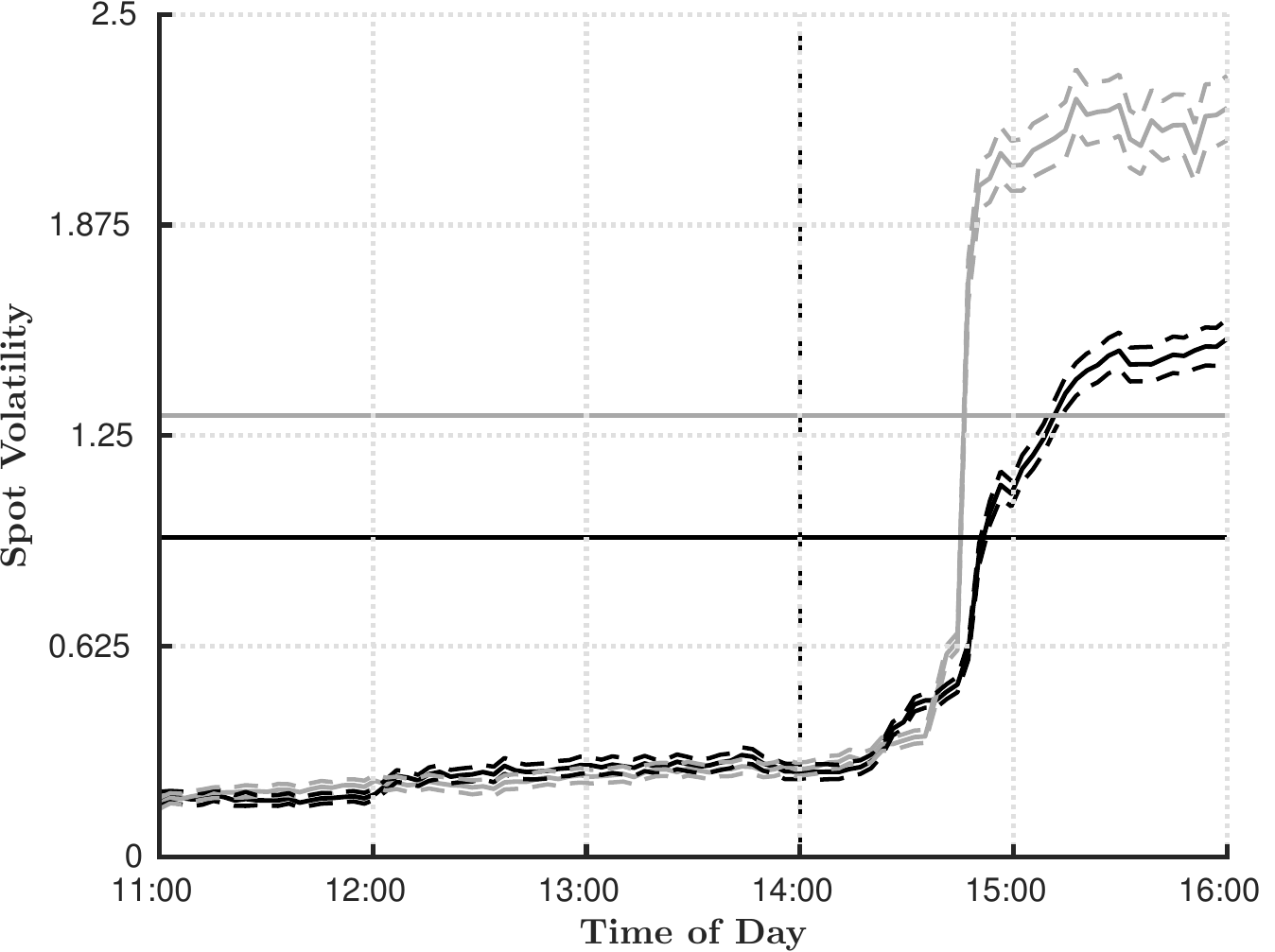}}\hfill
\caption{Spot covariances, correlations and volatilities  for AAPL and AMZN (05/06/10). In left plot, black lines (and left y-axis) represent correlations, grey lines (and right y-axis) covariances. In right plot, black lines are for AMZN, grey lines for AAPL. The dashed lines correspond to approximate pointwise $95\%$ confidence intervals according to Corollary~\ref{fclt}. The horizontal lines correspond to the cross-sectional median of \textit{integrated} covariance, correlation and volatility estimates. These are based on the LMM estimator of the integrated (open-to-close) covariance matrix by \citet{BHMR} accounting for serially dependent noise and using the same input parameter configuration as the spot estimators.  Covariances and volatilities are annualized.}
\label{fig:cov_corr_vol_AAPL_AMZN_20100605}
\end{figure}

The first study analyses the flash crash on May 6, 2010, see, e.g.~\citet{kiri2016}. Figure \ref{fig:cov_corr_vol_20100605} shows the cross-sectional deciles of spot covariances, correlations and volatilities on this day. We observe that correlations are virtually constant during the morning, but \emph{increase} slowly shortly after 2:00 pm, and, subsequently, \emph{decrease} quickly around 2:45 pm when prices began to return to their pre-crash levels. The latter is accompanied by an underlying pronounced \emph{increase} in covariances, while the deciles show that the cross-sectional distribution of covariances across all asset pairs is extremely skewed, revealing huge upward shifts in some covariances, but only very moderate reactions in others. Figure~\ref{fig:cov_corr_vol_20100605} also demonstrates that the corresponding reactions in spot volatilities have been much stronger, which explains the drop in median correlations from approx.~$0.5$ before 2:45 pm to approx.~$0.3$ right after 3:00 pm.  For comparison, Figure~\ref{fig:cov_corr_vol_AAPL_AMZN_20100605} displays the spot covariance, correlation and volatility estimates along with the corresponding approximate $95\%$ confidence intervals for AAPL and AMZN, which are the most liquid assets as measured by the number of mid-quote revisions. Most importantly, we observe that these particular spot covariance, correlation and volatility paths are in line with the patterns found in Figure~\ref{fig:cov_corr_vol_20100605}.  Accordingly, the latter are not an artifact of the cross-sectional aggregation across assets or asset pairs.

\begin{figure}
\centering
\hspace*{-4ex}\hfill\subfigure[Entire trading day]{\includegraphics[width=0.42\textwidth]{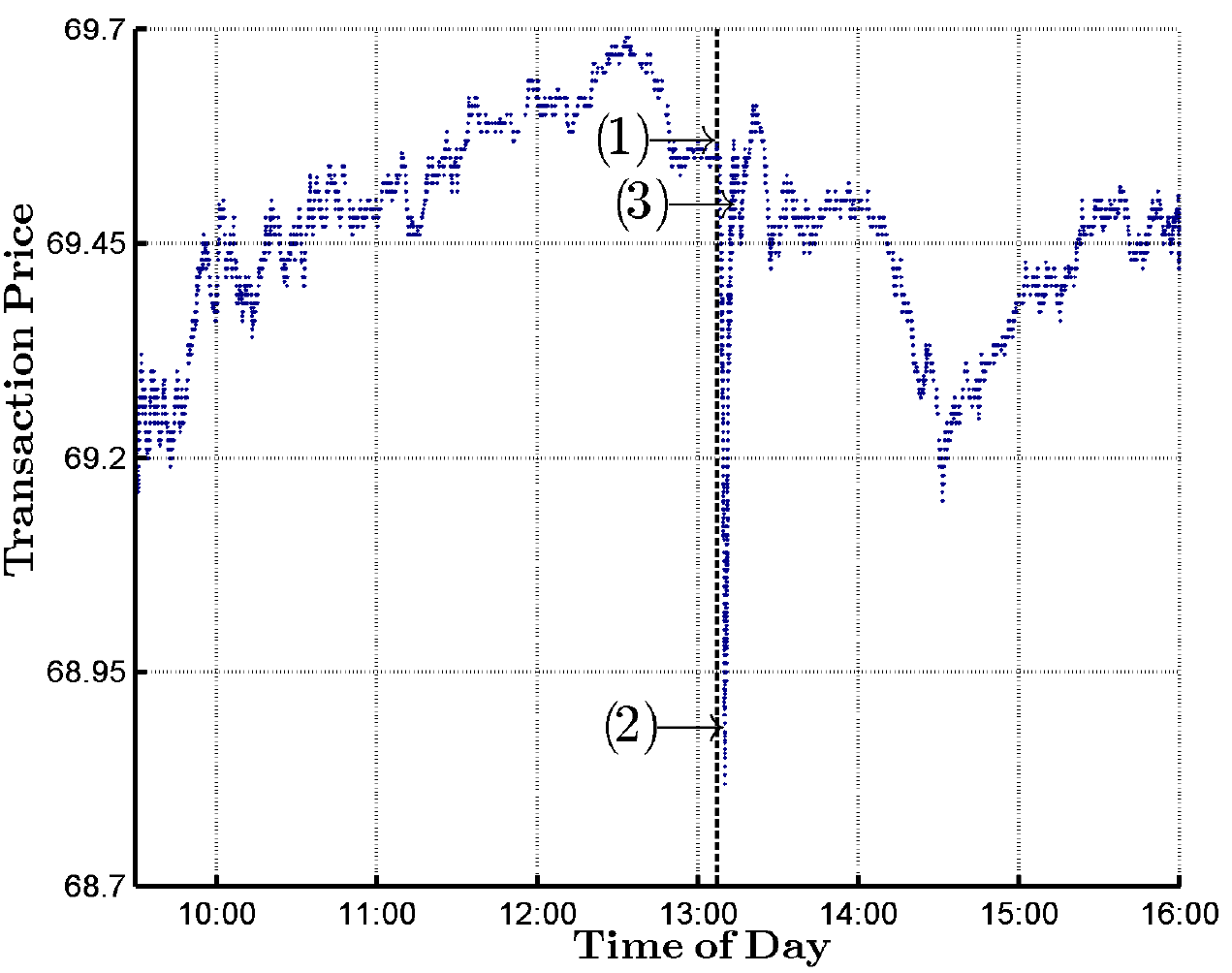}}\hspace*{-0.3ex}\hfill
\subfigure[1:00 pm  -- 1:30 pm]{\includegraphics[width=0.42\textwidth]{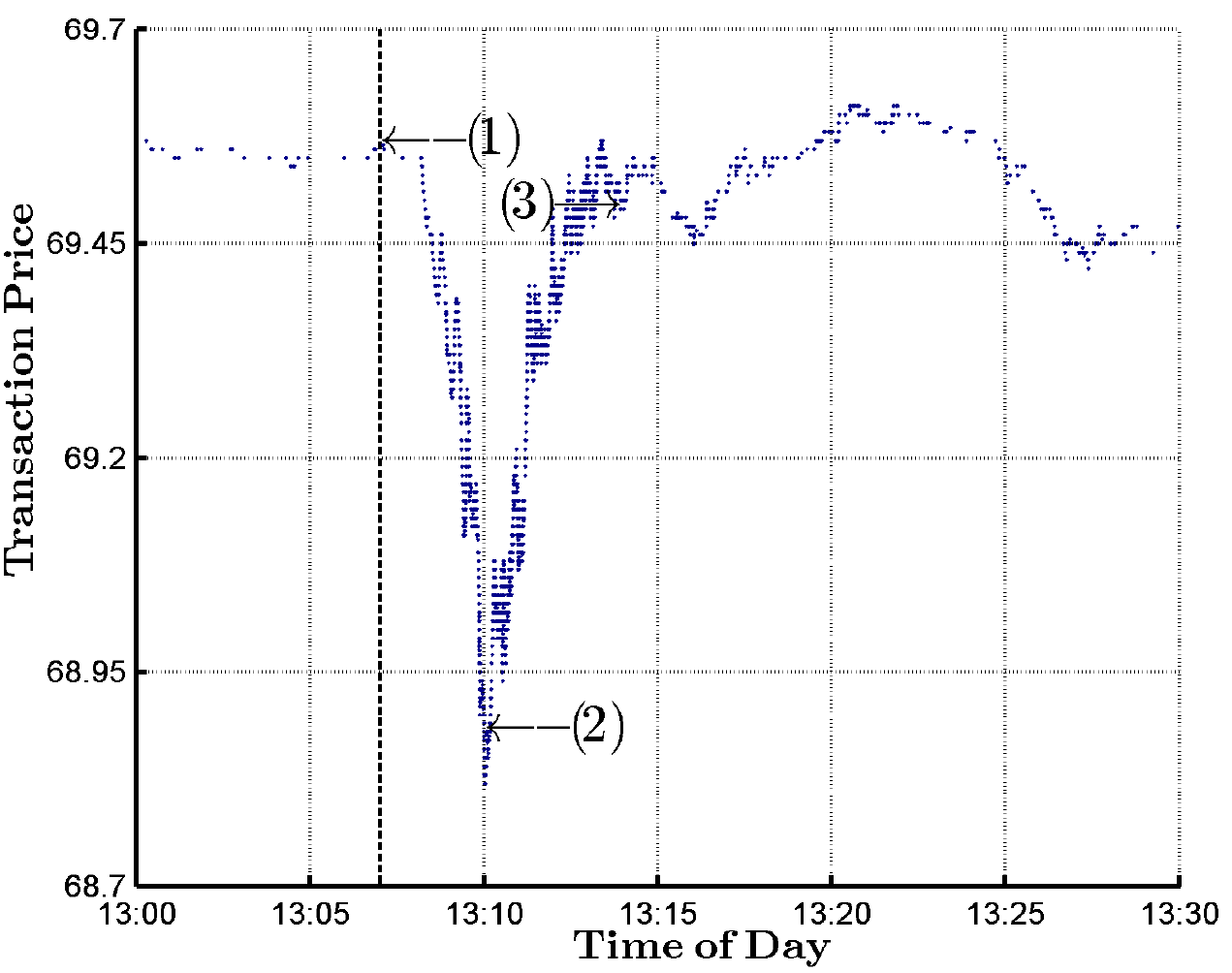}}\hfill
\caption{QQQ transaction prices (04/23/13). (1): Fake tweet from the account of AP stating ``Breaking: Two Explosions in the White House and Barack Obama is injured''. (2): Official denial by AP. (3): AP's twitter account suspended.}
\label{fig:price_20130423}
\end{figure}

\begin{figure}
\centering
\hspace*{-4ex}\hfill\subfigure[Covariances]{\includegraphics[width=0.42\textwidth]{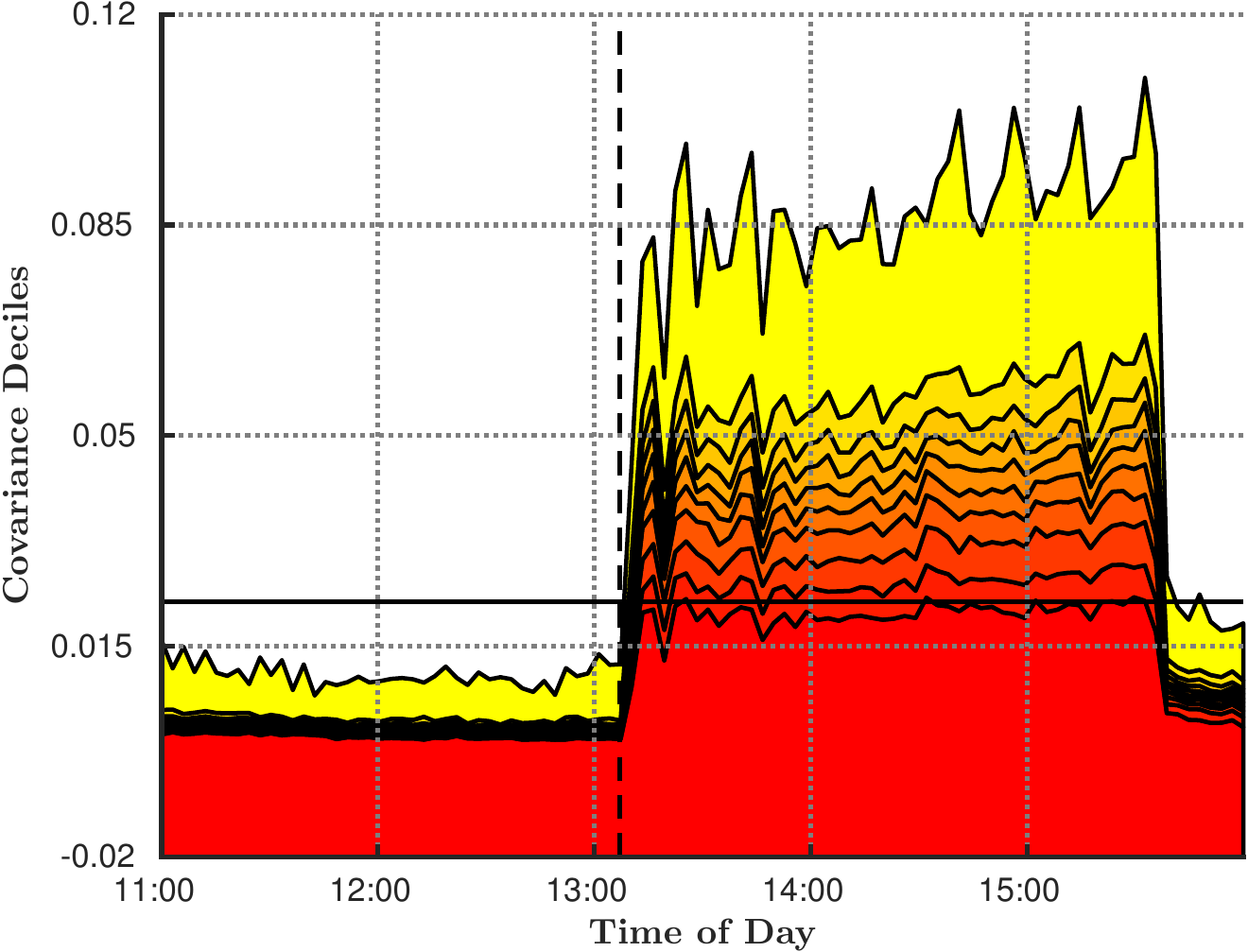}}\hspace*{-0.3ex}\hfill
\subfigure[Correlations]{\includegraphics[width=0.42\textwidth]{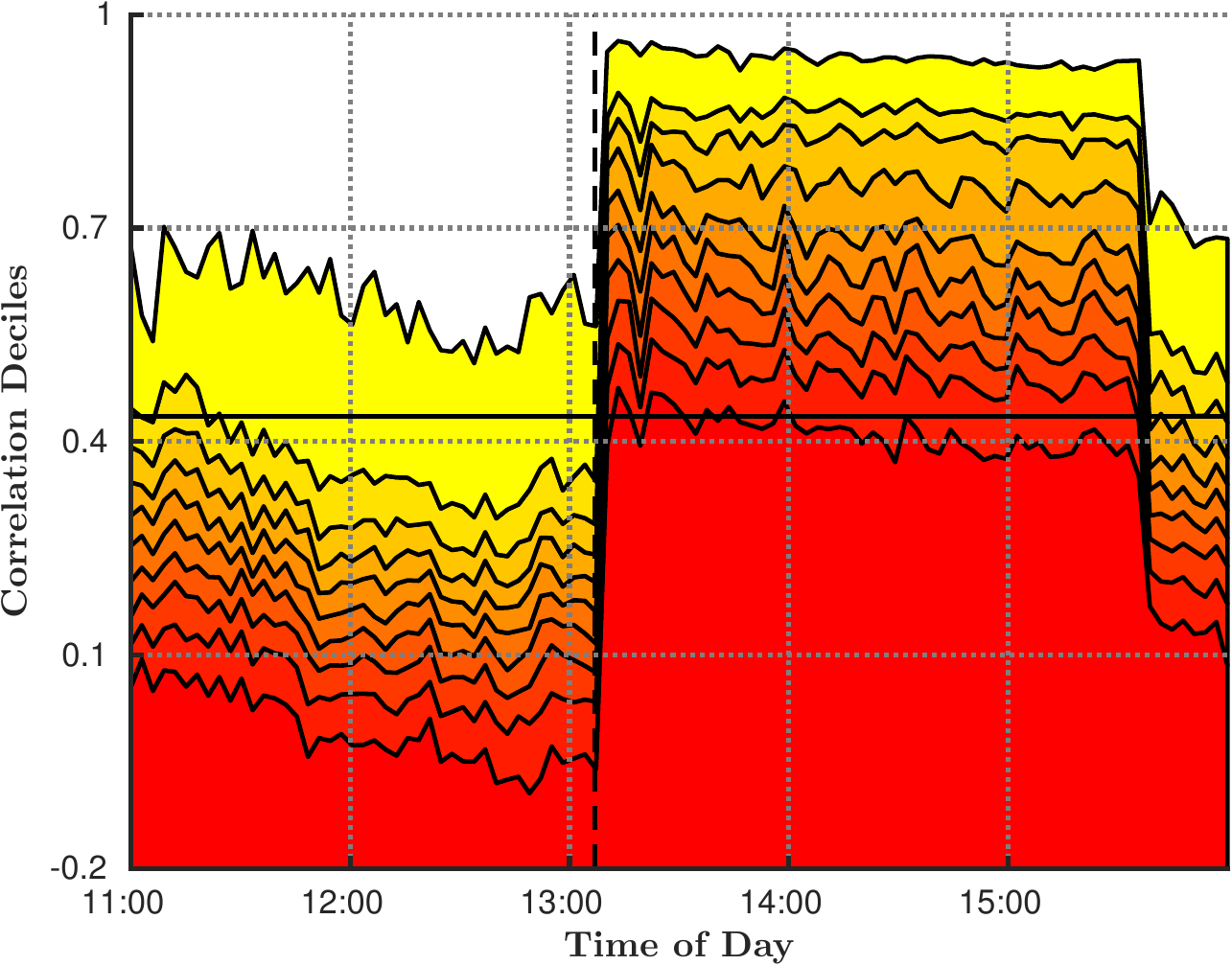}}\hfill\\
\subfigure[Volatilities]{\includegraphics[width=0.42\textwidth]{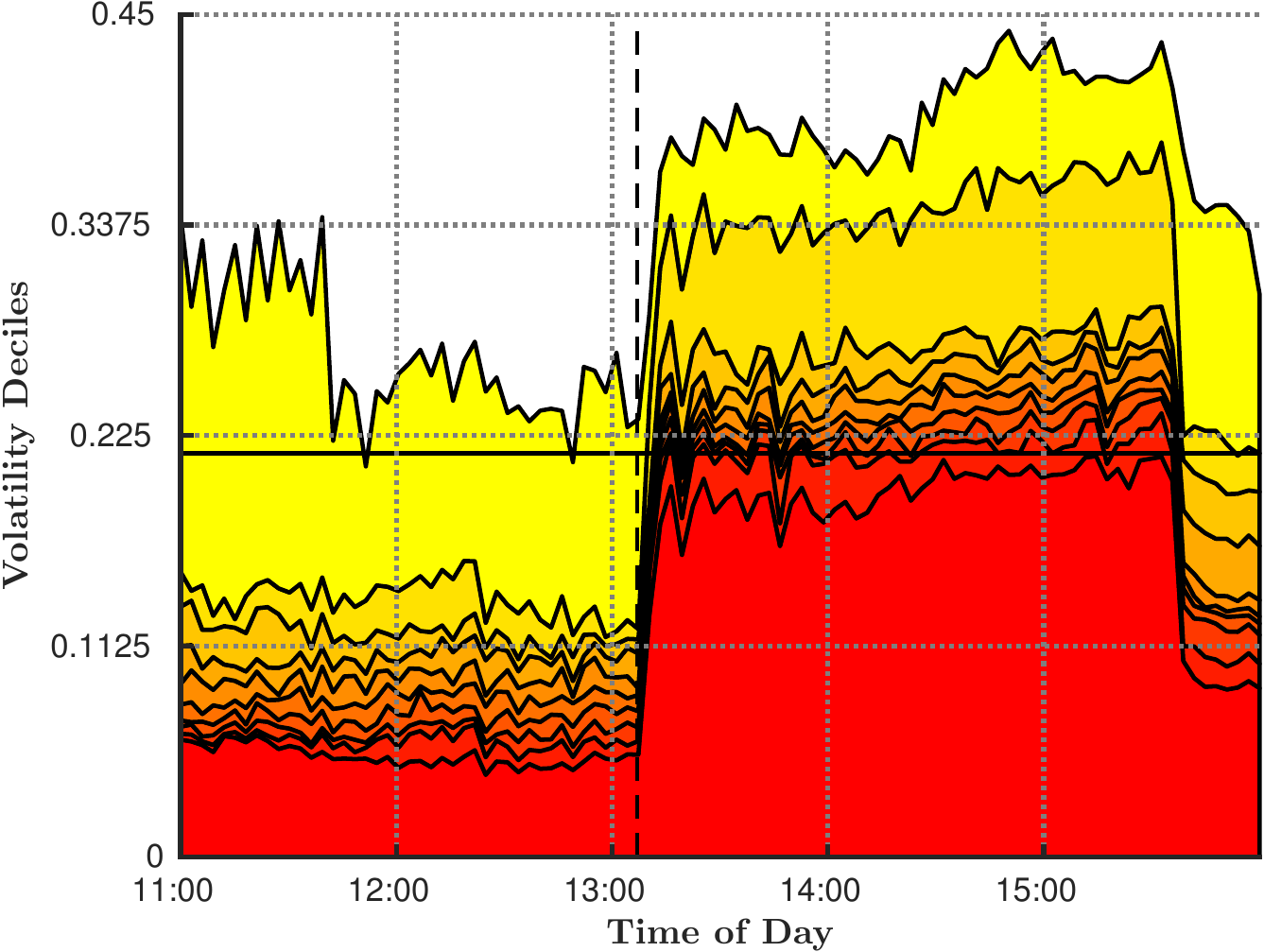}}\hfill
\caption{Cross-sectional deciles of spot covariances, correlations and volatilities (04/23/13).  The solid horizontal line corresponds to the cross-sectional median of \textit{integrated} covariance, correlation and volatility estimates. These are based on the LMM estimator of the integrated (open-to-close) covariance matrix by \citet{BHMR} accounting for serially dependent noise and using the same input parameter configuration as the spot estimators.  Covariances and volatilities are annualized.}
\label{fig:cov_corr_vol_20130423}
\end{figure}

\begin{figure}
\centering
\hspace*{-4ex}\hfill\subfigure[Covariances/Correlations]{\includegraphics[height=0.33\textwidth]{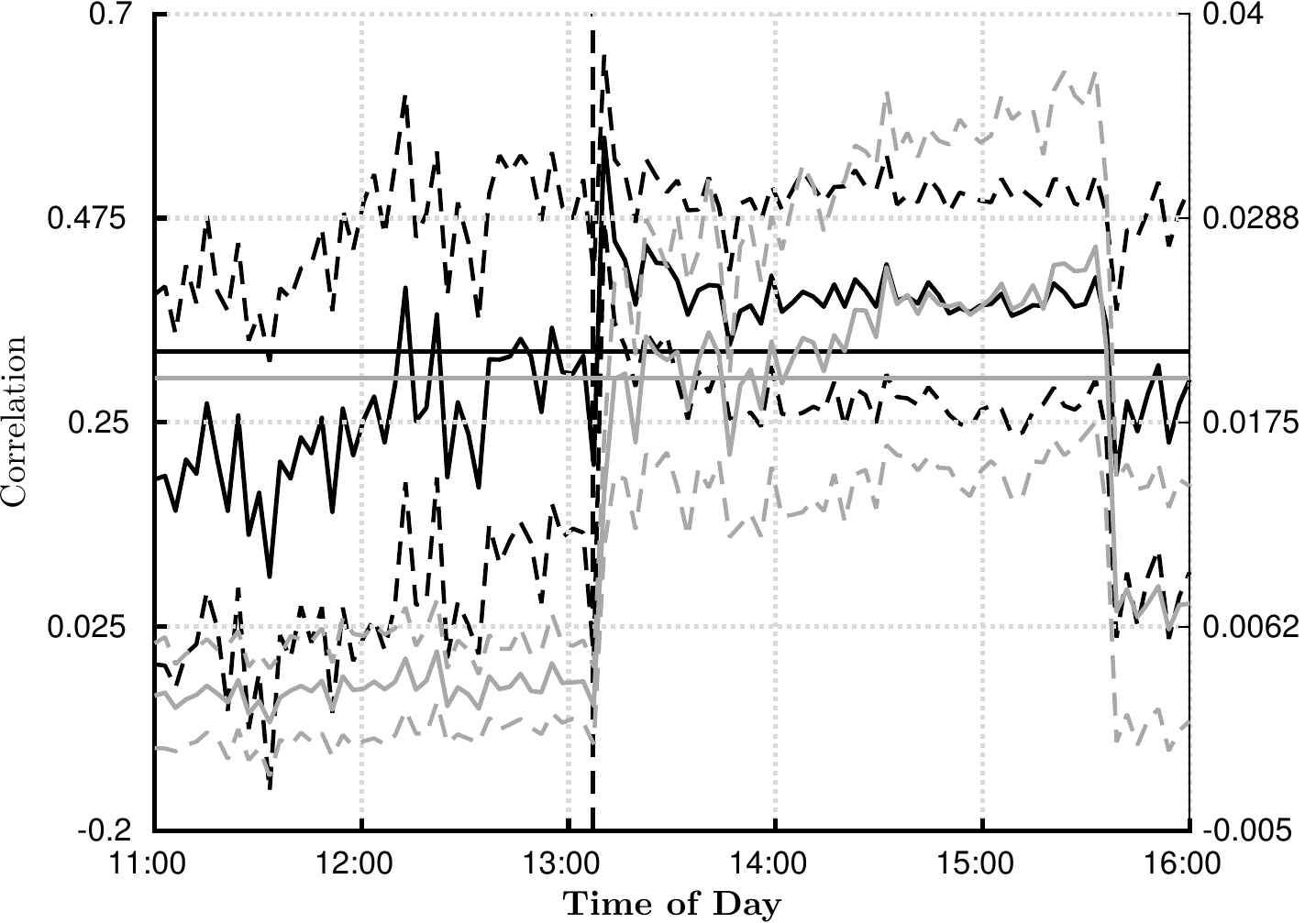}}\hfill
\subfigure[Volatilities]{\includegraphics[height=0.33\textwidth]{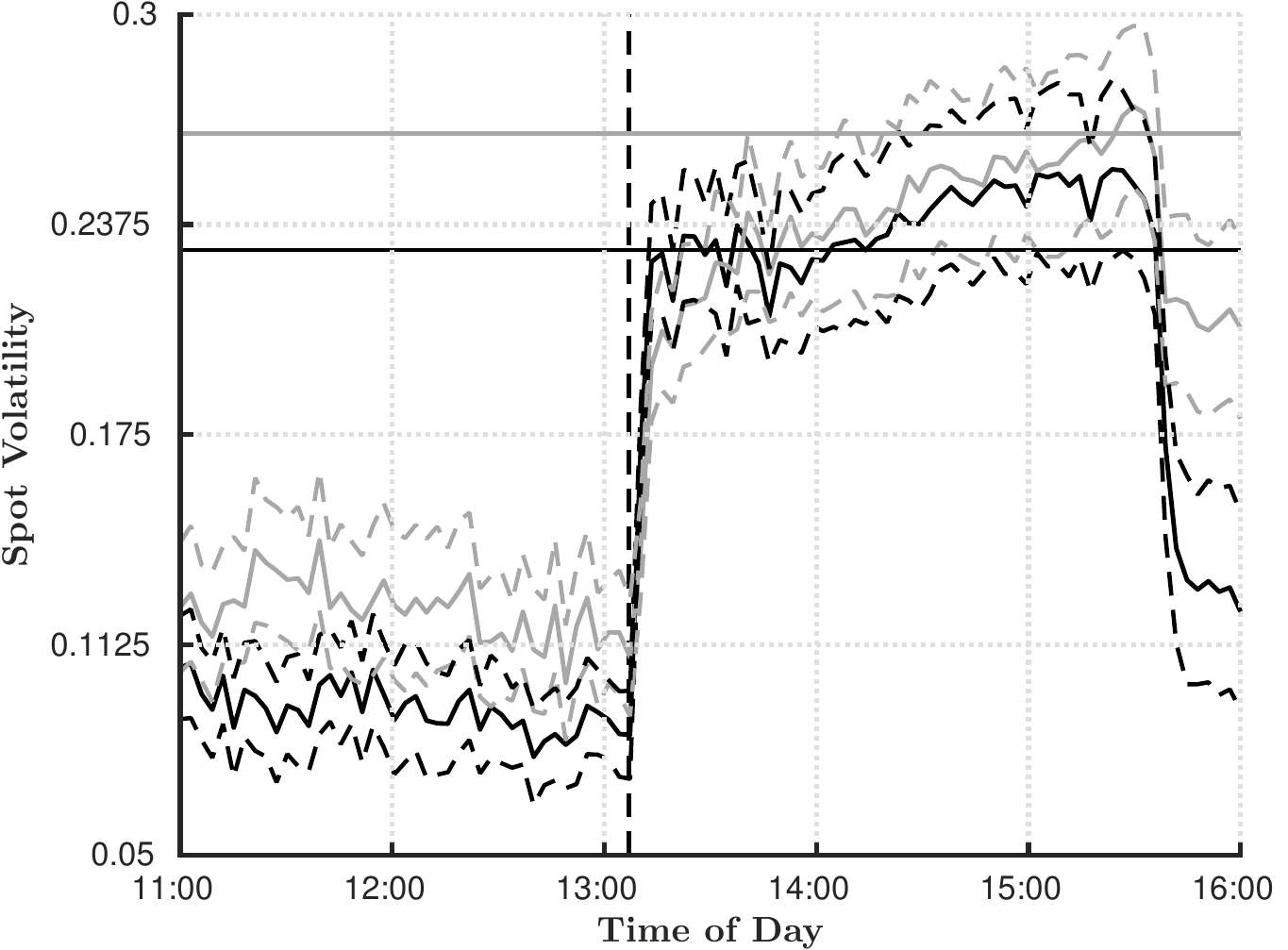}}\hfill
\caption{Spot covariances, correlations and volatilities  for AAPL and AMZN (04/23/13). In the left plot, black lines (and left y-axis) represent correlations, grey lines (and right y-axis) covariances. In the right plot, black lines are for AMZN, grey lines for AAPL. The dashed lines correspond to approximate pointwise $95\%$ confidence intervals according to Corollary~\ref{fclt}. The horizontal lines correspond to the cross-sectional median of \textit{integrated} covariance, correlation and volatility estimates. These are based on the LMM estimator of the integrated (open-to-close) covariance matrix by \citet{BHMR} accounting for serially dependent noise and using the same input parameter configuration as the spot estimators.  Covariances and volatilities are annualized.}
\label{fig:cov_corr_vol_AAPL_AMZN_20130423}
\end{figure}

Second, we analyze an  event, which is characterized by completely non-anticipated (and ultimately wrong) news. On 04/23/13 at around 1:07 pm, a fake tweet from the account of the Associated Press (AP) reported ``breaking'' news on two explosions in the White House, where the U.S. president (supposedly) got injured. At 1:10 pm, AP officially denied this message and  suspended its twitter account at 1:14 pm. Figure \ref{fig:price_20130423} shows the underlying price process and the timing of the corresponding events.  Our results in Figure \ref{fig:cov_corr_vol_20130423} show that (co-)variances and correlations strongly increase immediately after 1:07 pm.  The increase in covariances is stronger than for volatilities, whereby correlations (median) increase from approximately $0.2$ to $0.7$. The estimates suggest that the effect of elevated (co-)variances and correlations has been present for about two hours. As in the case of the May 2010 flash crash, this result is widely confirmed by the estimated paths for the specific asset pair of AAPL and AMZN displayed in Figure~\ref{fig:cov_corr_vol_AAPL_AMZN_20130423}. The above findings are remarkable given that the flash crash itself lasted only a couple of minutes and is similar to the effects observed during the May 2010 flash crash. Hence, effects of (flash) crashes on covariances may remain in the market for a considerable time period.

In summary, we conclude on the following findings: First, flash crash-type events cause abrupt upward movements in (co-)variances and correlations. 	Second, with prices ultimately returning to pre-shock levels, correlations move back, while the behavior of covariances and volatilities is more ambiguous. Depending on the nature of the market recovery process, they can both either decrease or increase. In the latter case,  the rise in volatilities is more pronounced, leading to reduced correlations. Our spot estimator seems to capture these effects quite well, since the observed reactions in the spot quantities are aligned with the timing of the underlying event. This indicates that the estimators are suitable to capture changes in dependence structures on a high time resolution.

\clearpage

\section{Conclusions\label{sec:concl}}
In this paper, we introduce an estimator for spot covariance matrices, which is constructed based on local averages of block-wise estimates of locally constant covariances. The proposed estimator builds on the local method of moments approach introduced by \citet{BHMR}. We show how to extend the LMM approach to the case of autocorrelations as well as endogeneities  in market microstructure noise and provide a suitable procedure for choosing the lag order in practice. For the resulting spot covariance matrix estimator, we derive a stable central limit theorem along with a feasible version that is straightforwardly applicable in empirical practice. An important result is that we are able to attain the optimal convergence rate, which is  $n^{1/8}$  under the assumption of a semi-martingale volatility matrix process with the efficient log-prices being subject to noise and a non-synchronous observation scheme.

Simulation exercises provide guidance on how to implement the estimator in practice and demonstrate its relative insensitivity with respect to the choice of block sizes, cut-offs and smoothing windows. Moreover, based on Nasdaq blue chip stocks, we provide detailed empirical evidence on the intraday behavior of spot covariances, correlations and volatilities.
In particular, we show that not only spot volatilities as previously documented in the literature, but also covariances and correlations reveal distinct intraday seasonality patterns. Further, we analyze how spot covariances change in periods of extreme market movements and show that intraday changes of (co-)volatility structures can be quite distinct and considerable. 

\section*{Acknowledgment and Web appendix}
\addcontentsline{toc}{section}{Acknowledgment and Web appendix}
In the web appendix, we provide the proofs, extended simulations and more detailed empirical results. Moreover, we provide commented code for the implementation of the methods.
We are grateful to the editor, an associate editor and two referees for helpful comments on a previous version.
\bibliography{ref}

\end{document}